\documentclass[12pt,a4paper,reqno]{amsart}

\usepackage{amssymb,mathrsfs,mathdots}
\usepackage{hyperref,verbatim}
\usepackage{colortbl,graphicx,overpic}
\usepackage{enumitem}
\usepackage{multirow}
\usepackage[margin=3cm]{geometry}

\newtheorem{theorem}{Theorem}[section]
\newtheorem{lemma}[theorem]{Lemma}

\newtheorem{proposition}[theorem]{Proposition}
{\theoremstyle{definition}

\newtheorem{remark}[theorem]{Remark}}
\newcommand{\Z}{\mathbb{Z}}
\newcommand{\R}{\mathbb{R}}
\newcommand{\CO}{\ensuremath{\mathcal{O}}}

\begin{document}

\title[Limit cycles of piecewise quadratic systems]{On the number of limit cycles in piecewise planar quadratic differential systems}

\author[F. Braun]{Francisco Braun}

\address{Departamento de Matem\'{a}tica, Universidade Federal de S\~ao Carlos, 
13565--905 S\~ao Carlos, S\~ao Paulo, Brazil}
\email{franciscobraun@ufscar.br}

\author[L.P.C. da Cruz]{Leonardo Pereira Costa da Cruz}

\address{Instituto de Ci\^encias Matem\'aticas e Computa{\c{c}}{\~a}o, Universidade de São Paulo,  13566--590 S\~ao Carlos, S\~ao Paulo, Brazil}
\email{leonardocruz@icmc.usp.br}

\author[J. Torregrosa]{Joan Torregrosa}

\address{Departament de Matem\`{a}tiques, Universitat Aut\`{o}noma de Barcelona, 08193 Be\-lla\-ter\-ra, Barcelona (Spain); Centre de Recerca Matem\`atica, Campus de Be\-lla\-ter\-ra, 08193 Be\-lla\-ter\-ra, Barcelona (Spain)}
\email{joan.torregrosa@uab.cat}

\subjclass[2010]{Primary: 34C07, 34C29, 34C25; Secondary: 37G15.}

\keywords{Periodic solution, averaging method, non-smooth differential system.}

\date{\today}

\begin{abstract}
We consider piecewise quadratic perturbations of centers of piecewise quadratic systems in two zones determined by a straight line through the origin. By means of expansions of the displacement map, we are able to find isolated zeros of it, without dealing with the unsurprising difficult integrals inherent in the usual averaging approach. We apply our technique to non-smooth perturbations of the four families of isochronous centers of the Loud family, $S_1$, $S_2$, $S_3$, and $S_4$, as well as to non-smooth perturbations of non-smooth centers given by putting different $S_i$'s in each zone. To show the coverage of our approach, we apply its first order, which recovers the averaging method of the first order, in perturbations of the already mentioned centers considering all the straight lines through the origin. 
Then we apply its second order to perturbations of the above centers for a specific oblique straight line. Here in order to argue we introduce certain blow-ups in the perturbative parameters. As a consequence of our study, we obtain examples of piecewise quadratic systems with at least $12$ limit cycles. By analyzing two previous works of the literature claiming much more limit cycles we found some mistakes in the calculations. Therefore, the best lower bound for the number of limit cycles of a piecewise quadratic system is up to now the $12$ limit cycles found in the present paper.
\end{abstract}

\maketitle

\section{Introduction} 

Consider the class of planar polynomial differential systems of degree $n$. The maximum number of \emph{limit cycles} that a system in this class can have is called the \emph{Hilbert number}, denoted by $H(n)$. The problem of finding $H(n)$ remounts to Hilbert and his $16$th problem which is, up to our knowledge, open until these days. Actually, although it is well known that $H(1) = 0$ and lower bounds for $H(n)$, $n\geq 2$, have been found over the years, upper bounds for it are still unknown for all $n\geq 2$. The search for $H(2)$ has been the object of intense study during the last century. The best-known lower bound for $H(2)$ was given by Shi \cite{Songling1980}, by means of an example of a quadratic differential system having $4$ limit cycles, so that $H(2)\geq4$. For the cubics, Li, Liu, and Yang \cite{LiLiuYang2009} showed that $H(3) \geq 13$. Denoting by $M(n)$ the maximum number of limit cycles bifurcating from a singular point of a polynomial system of degree $n$ as a degenerate Hopf bifurcation, it is clear that $M(n)$ is a lower bound for $H(n)$. Bautin \cite{Bautin1954} showed that $M(2) = 3$. Zoladek \cite{Zol1995,Zol2016} proved that $M(3) \geq 11$, see also a simpler proof by Christopher \cite{Chr2005}. In \cite{YuYun2014}, Yu and  Tian gave an example with $12$ limit cycles surrounding a singularity for cubic systems, so that $M(3)\geq 12$. This proof has some gaps but was corrected by Giné, Gouveia, and Torregrosa in \cite{GouTorreGine2021}.

We have witnessed an increasing interest in piecewise smooth systems. This is probably motivated by the wide range of applications that they have in modeling real phenomena, see, for instance, \cite{AcaBonBro2011,BerBudChaKow2008}. 
With a more theoretical point of view, it has been usual to explore the piecewise world by asking similar questions as in the smooth one \cite[Appendix A]{MR3838314}. 
This is the case of determining an analogous of the Hilbert number for piecewise polynomial systems, the heart of the present paper. In recent years some authors have obtained lower bounds for this new Hilbert number for low-degree systems. The aim of this paper is to study piecewise polynomial vector fields degree $2$ separated by a straight line. Before describing our results, for completeness and to situate the reader, we will detail the state of the art for degrees $1$ and $2$.

In this work, we consider the class of \emph{piecewise polynomial systems of degree $n$}: 
\begin{equation}\label{eq:0}
Z(x,y)= \left\{
\begin{array}{ccc}
Z^+(x,y),\quad h(x,y)\geq0,\\
Z^-(x,y),\quad h(x,y)\leq0,\\
\end{array}\right.
\end{equation}
where $Z^{\pm}$ are polynomial vector fields of degree $n$ and $h:\R^2 \rightarrow \R$ is the linear function 
$$
h(x,y) = ax+by, 
$$
$(a,b) \neq (0,0)$.  
So here the discontinuity curve $\Sigma = h^{-1}(0)$ is a straight line through the origin. 
Due to the fact that in $\Sigma$, the vector field $Z$ is bivaluated, we use the Filippov convention \cite{Fil1988} in order to define the local trajectories of $Z$ on $\Sigma$. 

In this context, we can consider a limit cycle of $Z$, i.e., an isolated periodic orbit of $Z$. A special one is the \emph{crossing limit cycle}, being a limit cycle intersecting $\Sigma$ trough the \emph{crossing region}. 

Here, analogously to the smooth case, we can consider the maximum number of limit cycles such a system can have. We denote this maximum number by $H_p(n)$. Precisely, $H_p(n)$ is the maximum number of limit cycles that a piecewise polynomial system of degree $n$ with a discontinuity curve being a straight line through the origin can have. Particularly we denote by $H_p^c(n)$ the maximum number of crossing limit cycles that a piecewise polynomial system in this class can have. Further, we denote by $M^c_p(n)$ the maximum number of crossing limit cycles bifurcating from a singular point or sliding set that a system in this class can have. 

Up to now, for piecewise polynomial systems of degree $1$, an example with $3$ limit cycles was first detected numerically by Huan and Yang \cite{HuanYang2010}. 
Soon later it was analytically proved by Llibre and Ponce \cite{LliPon2012}. On the other hand, examples of piecewise systems of degree $1$ having more than $3$ crossing limit cycles are not known. The existence of $3$ crossing limit cycles in piecewise systems of degree $1$ was also obtained using different techniques, among others, by means of perturbations of centers, by degenerated Hopf bifurcations from infinity, or as an application of Poincar\'e--Miranda Theorem. See more details in \cite{BuzPesTor2013,FrePonTorTor2021,FrePonTor2014,GasMan2020,LliNovTei2015a}. So $H^c_p(1) \geq 3$.  

For $n=2$, Llibre and Mereu \cite{MereuLlibre2014} obtained at least $5$ limit cycles by perturbing a suitable quadratic isochronous center and applying averaging theory of first order. Llibre and Tang \cite{LiTan2019} obtained $8$ crossing limit cycles by using averaging theory of order $5$ in a quadratic perturbation of the linear center. Tian and Yu \cite{YunPei2015} claimed the existence of $10$ limit cycles in a quadratic piecewise system. Actually, with the customary additional limit cycle coming from the pseudo-Hopf bifurcation, see Section~\ref{pseudo}, we can assert that these authors indeed found $11$ limit cycles. 

Then by perturbing a suitable isochronous quadratic center, the birth of $16$ crossing limit cycles was claimed in \cite{CruzNoaesTorre2019}. We redo the calculations of \cite[Theorem 1.1]{CruzNoaesTorre2019}, with the same technique used there, and could not reply this number of limit cycles, so we confirm that some mistake occurred in the application of averaging theory of order $2$ there. Therefore, we can assure that such result is not correct.  Also, in \cite[Proposition 3.1]{GouTorre2021}, it was claimed the appearance of $13$ crossing limit cycles after a non-smooth quadratic perturbation of a suitable quadratic center. We redo the calculations, with the same technique used there, and confirm that some mistakes occurred as well. 

In this paper, we analyze non-smooth quadratic perturbations of suitable quadratic smooth and non-smooth centers. As a consequence of our study, we provide a new lower bound for $H^c_p(2)$: $12$ limit cycles. This is up to our knowledge the best lower bound for $H^c_p(2)$ up to now. 

The quadratic smooth centers we perturb are the isochronous quadratic systems $S_1$, $S_2$, $S_3$ and $S_4$ of the Loud \cite{Lou64} family, written after suitable linear changes of variables, as  
\begin{equation}\label{eq:10}
\begin{aligned}
& S_1:
\left\{\begin{array}{l}
\dot x= -y+x^2-y^2,\\[4pt]
\dot y=x(1+2y).
\end{array}
\right. 
& \phantom{4569} & S_2: 
\left\{\begin{array}{l}
\dot x=-y+x^2,\\[4pt]
\dot y=x(1+y).
\end{array}
\right. 
\\[4pt]
& S_3: 
\left\{\begin{array}{l}
\dot{x} = -y-\frac{4}{3}x^2,\\[4pt]
\dot{y} = x \big(1-\frac{16}{3}y\big).
\end{array}
\right. 
& \phantom{4569} & S_4: 
\left\{\begin{array}{l}
\dot x=-y+\frac{16}{3}\,x^2-\frac{4}{3}\,y^2,\\[4pt]
\dot y=x(1 + \frac{8}{3} y). 
\end{array}
\right.
\end{aligned}
\end{equation} 
The forms presented here are the ones of \cite{ChaSab1999}. 
See also \cite{MarRouTon1995}. 

Perturbations of the centers $S_1$, $S_2$, and $S_3$ frequently appear in the literature. For instance, the first order averaging method was used in the above-mentioned paper \cite{MereuLlibre2014} to obtain $4$ and $5$ limit cycles by perturbing $S_1$ and $S_2$, respectively, inside the class of piecewise quadratic systems with two zones separated by the straight line $y=0$. 

But since the usual way to address the problem with averaging theory is to consider suitable linearizations of the center to be perturbed, there are few approaches by using different lines than the horizontal or vertical ones, as well as perturbations of $S_4$. The reason is that the known linearizations of $S_1$, $S_2,$ and $S_3$ do not send all the straight lines onto straight lines, as well as there are no birational linearizations of $S_4$ and so it is very difficult to deal with the problem. Our approach does not use linearizations, so we can apply it to any straight line as well as to $S_4$. 

Before stating our results, we fix some notation. 
From now on, we assume that the origin of coordinates is a non-degenerate center of a piecewise polynomial vector field $Z$ of degree $n$ as in \eqref{eq:0}. 
A \emph{non-smooth polynomial perturbation of degree $n$} of the center $Z$ is the vector field
\begin{equation}\label{eq:8}
Z_\varepsilon (x,y)=\left\{
\begin{array}{l}
\hspace{-.15cm} Z^+(x,y) + \varepsilon Z_1^+(x,y), \textrm{ if } h(x,y) \geq 0,\\[4pt]
\hspace{-.15cm} Z^-(x,y) +  \varepsilon Z_1^-(x,y), \textrm{ if } h(x,y) \leq 0,
\end{array}
\right. 
\end{equation}
where $\varepsilon > 0$ and $Z_1^\pm(x,y) = \left(P_1^\pm(x,y),Q_1^\pm(x,y)\right)$ are polynomial vector fields of degree $n$ without constant terms. As it will be clear below, with the technique used in this paper, we do not need the usual perturbations with higher orders in $\varepsilon$. That is, we do not improve the number of limit cycles by adding such perturbations. 

In order to fix the notation from now on, we write 
\begin{equation}\label{perturbation_piecewise}
P_1^{\pm}(x,y) = \sum_{0<|\eta|\leq n} a^\pm_{\eta} x^{\eta_1} y^{\eta_2},\quad Q_1^{\pm}(x,y) = \sum_{0<|\eta|\leq n} b^\pm_{\eta} x^{\eta_1} y^{\eta_2}, 
\end{equation}
where $\eta = (\eta_1, \eta_2)$ is a pair of non-negative integer numbers and $|\eta| = \eta_1 + \eta_2$. 

We remark that we do not consider constant terms in the perturbations as it is known that at least one limit cycle is born by adding suitable constant terms a posteriori by means of a \emph{pseudo-Hopf type bifurcation}, as recalled in Section~\ref{pseudo} below. 

We consider the rational parametrization of $\mathbb{S}^1\cap\{x\geq 0\}$ given by 
\begin{equation}\label{parameter}
\tau\mapsto \left(\frac{1-\tau^2}{1+\tau^2}, \frac{2 \tau}{1+\tau^2} \right), 
\end{equation}
$\tau \in [-1,1]$. 
So the angle $\alpha \in [-\pi/2, \pi/2)$ can be changed by this parametrization varying $\tau$ in $[-1,1)$. 
The straight line $h(x,y) = a x + b y = 0$, $b\geq 0$, with slope $\alpha$, i.e. $a/\sqrt{a^2+b^2} = - \sin \alpha$ and $b/\sqrt{a^2+b^2} = \cos \alpha$, can be identified with $\tau\in [-1,1)$ by the equations $\sin \alpha = \frac{2 \tau}{1+\tau^2}$ and $\cos \alpha =  \frac{1-\tau^2}{1+\tau^2}$. 

In Theorem~\ref{theorem:order1} we apply the order $1$ of our method to quadratic perturbations of $Z= S_1$, $S_2$, $S_3,$ and $S_4$ considering all the values of $\tau\in [-1,1)$, i.e. all the straight lines through the origin. This is equivalent to the usual averaging theory of first order. For higher orders, our approach is quite different, as we explain in a while.

\begin{theorem}\label{theorem:order1}
Considering the straight line $\tau\in (-1,1)\setminus\{0\}$, the number of crossing limit cycles bifurcating from the origin after applying averaging theory of order $1$ on non-smooth quadratic perturbations of the isochronous $S_1$, $S_2$, $S_3$, and $S_4$ is at least $7$, $8$, $9$, and $9$, respectively, while for $\tau=0$ or $\tau=-1$, the numbers are $5$, $6$, $7$, and $6$, or $5$, $6$, $5$, and $6$, respectively. 
\end{theorem}

In the next theorem, we also consider quadratic perturbations for all values of $\tau$, but now for $Z$ such that $Z^+ = S_1$ and $Z^- = S_2$. This is possible because $Z$ defined in this manner has a center at the origin. Indeed, we will explain this latter, in Lemma~\ref{lemma:center}. We emphasize that, as far as we know, this is the first time a non-smooth center is studied by a perturbation. 

\begin{theorem}\label{theorem:main3}
Considering the straight line $\tau \in (-1,1)\setminus\{0\}$, the number of crossing limit cycles bifurcating from the origin after applying averaging theory of order $1$ on non-smooth quadratic perturbations of the system $Z$ such that $Z^+ = S_1$ and $Z^- = S_2$ is at least $10$, while for $\tau = 0$ or $\tau=1$, the numbers reduce to $6$ or $8$, respectively. 
\end{theorem}

Now, to simplify computations, we fix $\tau = 1/2$ in order to consider the straight line with slope $\alpha$ such that $\left(\cos \alpha, \sin \alpha \right)$ is the point with rational coordinates $\left(3/5, 4/5 \right)$ and use our method with order $2$, obtaining: 

\begin{theorem}\label{theorem:Main2}
Considering $\tau=1/2$, after non-smooth perturbations of order $2$ of $S_1$, $S_2$, $S_3$, and $S_4$, we obtain at least $7$, $10$, $9$, and $12$ crossing limit cycles bifurcating from the origin, respectively. 
\end{theorem} 

The same result can be obtained for other values of $\tau$, but we do not get it for any value of $\tau$ because of high computational costs. 

\begin{theorem}\label{theorem:main} 
Considering $\tau=1/2$, after non-smooth perturbations of order $2$ of the piecewise quadratic system given by $Z^+ = S_1$ and $Z^- = S_2$ we obtain at least $12$ crossing limit cycles bifurcating from the origin. 
\end{theorem} 
The following table summarizes our lower bounds for the number of limit cycles bifurcating from the origin up to order $2$, with $\tau=1/2$. 

\begin{center}
\begin{tabular}{|c|c|c|}
\cline{2-3}
\multicolumn{1}{c|}{} & $1$st order & $2$nd order  \\
\hline  
$S_1$ & $7$ & $7$  \\
\hline 
$S_2$ & $8$ & $10$  \\
\hline 
$S_3$ & $9$ & $9$  \\
\hline 
$S_4$ & $9$ & $12$  \\
\hline
$S_1\&S_2$ & $10$ & $12$  \\
\hline
\end{tabular}
\end{center}

This paper provides lower bounds for the Hilbert number by means of non-smooth perturbations of some quadratic smooth and non-smooth centers. As it is widely known, the main difficulty in this kind of study, either in the smooth or in the non-smooth case, is the huge amount of calculations one has to deal with. By using the usual averaging techniques, it is always a big challenge to calculate the integrals appearing in the procedure. Anyway, after obtaining the averaged functions, one has to analyze isolated zeros of it. We make a remark towards explaining how to achieve the results of this paper with less effort than one is used to see in similar works. 

When one is interested only in zeros of the averaged function close to $0$, which is the case when analyzing the birth of limit cycles after perturbing a center around an equilibrium point, it is common to use Taylor series expansion of the averaged function around $0$ in order to reason. It is precisely here we propose a different but almost equivalent approach: we expand in Taylor series in one of the variables \emph{before} calculating the integrals in the other variable. So invariably we obtain a sum of very simple integrals. It is worth noting that the idea of this work is developed for concrete cases in two dimensions, however, it is possible to develop it in more generic cases in dimensions greater than or equal to three. The issue we face now is that we cannot always assure that the first averaged function is identically zero for applying the second order as it is usual. This is clear because we are only calculating the expansion in the Taylor series of the first order function, and if a Taylor series has its first $k$ terms equal to zero, this of course does not imply the whole series is identically zero. When this happens, i.e. the first order is not identically zero, for obtaining results using second order functions, we have to deal with a combination of first and second terms, what we actually do. Of course, the same is true if we want results by using yet higher orders. So our method agrees with the usual averaging method for the first order. And it only agrees with the $n$th method of averaging if we are able to assure that the less order functions are identically zero. 

Not only our aim in this paper is to call attention to a correction in the up to now known lower bound for $H_p^c(2)$, but also we present the above mentioned alternative way of reasoning. We then apply our method to obtain the above stated theorems, illustrating their coverage. We observe that up to our knowledge, there are no results upon perturbating $S_4$, mainly because there does not exist a birational linearization of the center at the origin, making the integral calculations prohibitive, a technical tool we can now dismiss, and we include perturbations of $S_4$. 

We organize the content of this paper as follows: in Section~\ref{section:averaging} we recall the results for finding limit cycles and explain our expansion method in detail. In Section~\ref{section:S1toS4} we apply our technique to obtain results upon non-smooth perturbations of $S_1$, $S_2$, $S_3$, and $S_4$, as well as the combination of $S_1$ and $S_2$, up to order $1$ for all the straight lines crossing the origin, proving Theorems~\ref{theorem:order1} and \ref{theorem:main3}. Finally, in Section~\ref{section4}, we apply the second order of our method to perturbations of the same centers but now only with the straight line $\tau=1/2$: here we prove Theorems~\ref{theorem:Main2} and \ref{theorem:main}. It is in the proofs of these theorems that we explain how to deal with the mentioned combination of orders $1$ and $2$. 

We finish this introduction section remarking that up to our knowledge, this is the first time a ``pure-discontinuous'' center is perturbed and studied, in Theorems~ \ref{theorem:main3} and \ref{theorem:main}. 

\section{Finding limit cycles of piecewise polynomial planar systems by perturbing a piecewise center}\label{section:averaging}

\subsection{The difference function}\label{subsection:averaging}
We begin by applying polar coordinates $(x,y)=(r\cos\theta,r\sin\theta)$, so that the perturbed center \eqref{eq:8} writes
$$
\widetilde Z_{\varepsilon}(\theta,r)=\left\{
\begin{array}{l}
\!\!\widetilde Z^+_\varepsilon(\theta,r),\quad \textrm{if}\quad  \alpha\leq\theta\leq\alpha+\pi,\\
\!\!\widetilde Z^-_\varepsilon(\theta,r), \quad \textrm{if}\quad \alpha-\pi\leq\theta\leq\alpha. 
\end{array}
\right. 
$$ 
Then we make $\theta$ as the new independent variable, and the differential equation associated with this vector field becomes 
\begin{equation}\label{eq:drdt}
r'(\theta)=\frac{dr}{d\theta}=F_0(\theta,r) + \sum_{i=1}^m \varepsilon^i F_i(\theta,r) + \CO(\varepsilon^{m+1}),
\end{equation}
where
\[
F_i(\theta,r)=\left\{
\begin{array}{ccc}
\!\! F_i^+(\theta,r),&\textrm{if}& \alpha\leq\theta \leq\alpha+\pi,\\
\!\! F_i^-(\theta,r),&\textrm{if}&\alpha-\pi\leq\theta\leq\alpha,
\end{array}\right.
\]
being $F_i^{\pm}:[\alpha-\pi,\alpha+\pi]\times (0,\rho^*)\rightarrow\R$ analytic, for small enough $\rho^*$ (because $Z$ is a center) and $2\pi$-periodic in the variable $\theta$ for $i=0,1,\ldots, m$. By means of the change $\theta \mapsto \theta +\alpha$, we can assume $\alpha = 0$. 

We denote by $\varphi_{\varepsilon}(\theta,r)$ the flow of \eqref{eq:drdt} and write 
$$
\varphi_{\varepsilon}(\theta,r) = \sum_{i=0}^m \varepsilon^i \varphi_i(\theta,r) + \CO(\varepsilon^{m+1}), 
$$
so that 
$$
\sum_{i = 0}^k\varepsilon^i \varphi_i(\theta,r)
$$ 
is the $k$-jet of the solution of the initial value problem
$$ 
z'(s) = \sum_{i = 0}^k \varepsilon^i F_i(s, z(s)) + \CO(\varepsilon^{k + 1}),\quad z(0) = r,  
$$
$k = 0,\ldots, m$. 
Each function $\varphi_i: \R \times [0,\infty) \to [0,\infty)$ is of the form $\varphi_i^+(\theta, r)$ for $\theta \in [2 \kappa \pi, (2\kappa+1)\pi]$, and $\varphi_i^-(\theta, r)$ for $\theta \in [(2 \kappa - 1)\pi, 2 \kappa \pi]$, $\forall \kappa \in \Z$, $\theta$ eventually restricted to maximal intervals of solutions. 

Beginning with $\varphi_0^{\pm}$, the solution of the initial value problem 
\begin{equation*}
z'(s) = F_0^{\pm}(s,z(s)),\quad z(0) = r, 
\end{equation*}
we get that $\varphi_i^{\pm}$, $i=1,\ldots, k$, are given by the integral equations, according to, for instance, an adaptation of \cite{MR3177572} for the non-smooth case, 
\begin{equation}\label{y1y2y3}
\begin{aligned}
\varphi_1^\pm(\theta,r) &= \int_0^\theta \left(F_1^\pm(s,\varphi^\pm_0(s,r)) + \frac{\partial F^\pm_0}{\partial r}(s,\varphi^\pm_0(s,r)) \varphi_1^{\pm}(s,r) \right) ds, \\
\varphi_2^\pm(\theta,r) & = \frac{1}{2} \int_0^\theta\left(2 F_2^\pm(s,\varphi^\pm_0(s,r)) + 2 \frac{\partial F_1^\pm}{\partial r}(s, \varphi^\pm_0(s,r)) \varphi_1^\pm(s,r) \right. \\ 
& \phantom{={}} \left. + \frac{\partial^2 F_0^\pm}{\partial r^2}(s,\varphi^\pm_0(s,r)) \varphi_1^{\pm}(s,r)^2 + \frac{\partial F^\pm_0}{\partial r}(s,\varphi^\pm_0(s,r)) \varphi_2^{\pm}(s,r) \right)ds, 
\end{aligned}
\end{equation}
and $\varphi_i^{\pm}$'s, for $i\geq3$, are given recursively adapting \cite{MR3177572}. 

We define the \emph{$i$-difference function} as 
\begin{equation*}\label{media}
\delta_i(r) = \varphi_i^+(\pi,r) - \varphi_i^-(-\pi,r), 
\end{equation*}
$i=1,\ldots, k$, so that 
$$
\sum_{i=1}^k \varepsilon^i \delta_i(r)
$$ 
is the $k$-jet of the \emph{difference function} 
$$
\Delta(r,\varepsilon) = \varphi_{\varepsilon}(\pi,r) - \varphi_{\varepsilon}(-\pi,r). 
$$ 
Evidently, as the origin of the non-perturbed vector field is a center, $\varphi_0^+(\pi,r) - \varphi_0^-(-\pi,r) = 0$, so $\delta_0$ does not appear. 

We remark that the first non-vanishing coefficient of the Taylor expansion of the difference function provides the stability as it does the equivalent one for the usual displacement function, but with opposite signs. 

Clearly, each simple zero of $r\mapsto \Delta(r, \varepsilon)$ provides a hyperbolic $2 \pi$ periodic solution of \eqref{eq:drdt}, and so a hyperbolic limit cycle of \eqref{eq:8}. By the implicit function theorem, for $\varepsilon$ small enough, for each simple zero $\overline{r}$ of $\delta_{k}(r)$, where $\delta_i(r)\equiv 0$, $i = 1,\ldots, k - 1$, there exists a simple zero of $r \mapsto \Delta(r, \varepsilon)$ converging to $\overline{r}$. In this situation, looking for simple zeros of $\delta_{k}(r)$ in order to determine limit cycles of \eqref{eq:8} is usually termed as \emph{Averaging Method of Order $k$}, see \cite{MR3177572}, for instance. 

After a glance at the formulas of $\varphi_i$ above, it certainly comes with no surprise to the reader that finding the explicit expression of the $i$-difference functions $\delta_i(r)$ is a challenging task. In order to simplify this, it is usual to consider systems where $F_0 \equiv 0$. In the case of isochronous non-degenerate centers, it is always possible to find an analytic linearization that transforms it into the linear center, in which case $F_0 \equiv 0$. Even so, the integrals one has to calculate are difficult. Moreover, manageable linearizations do not always exist. For instance, the quadratic isochronous centers $S_1$, $S_2$, and $S_3$ have birrational linearizations, but $S_4$ has not, see \cite{ChaSab1999}. 

We propose a milder approach: instead of looking for explicit formulas for $\delta_i(r)$, we look for Taylor expansions in $r$. As we are going to see right below, in Section~\ref{subsection_aproximation}, this produces an algorithm where we have just to integrate trigonometric functions. But since we will only have a jet of $\delta_i(r)$, this may be a problem when we want to deal with order $k > 1$, as we will never be sure that $\delta_i(r) \equiv 0$, $i = 1,\ldots, k-1$, before looking for zeros of $\delta_k(r)$. To deal with this, we propose a \emph{blowing up} technique we will explain in detail in the proofs of Theorems~\ref{theorem:Main2} and \ref{theorem:main}.

Before explaining the algorithm of the mentioned expansion, we recall the so-called \emph{pseudo-Hopf type bifurcation}, that provides at least one more limit cycle in the non-smooth case by adding constant terms in the perturbation \eqref{eq:8}. 

\subsection{Pseudo-Hopf type Bifurcation}\label{pseudo}
\begin{figure}[!h]
\begin{center}
\begin{overpic}[scale=.65]{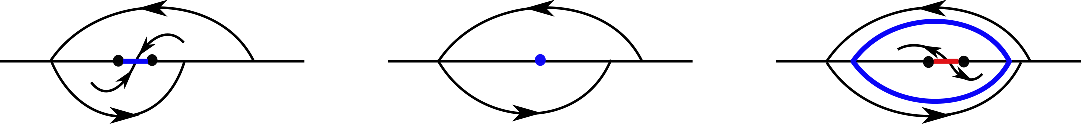}
\put(7,-3){$b<0$}
\put(45,-3){$b=0$}
\put(82,-3){$b>0$}
\end{overpic}
\end{center}
\caption{Pseudo-Hopf type bifurcation.}\label{fi:pseudohopf}
\end{figure}

Roughly speaking, in a Hopf type bifurcation, see for instance \cite{HalKoc1991}, a limit cycle arises when a monodromic equilibrium point changes its stability. 
In piecewise differential systems, a pseudo-Hopf type bifurcation describes an analogous phenomenon, but now a limit cycle arises when the sliding segment changes its stability. In this case, the size of the sliding segment takes the role of the trace in the analytic context. This phenomenon was described firstly in \cite{Fil1988}, but called pseudo-Hopf in \cite{Kuznetsov2003} when the pseudo-equilibrium point is of fold-fold type with codimension $1$. For codimension $2$, see \cite{GuarSeaTei2011}. The proof, as in the classical Hopf bifurcation, is a direct consequence of an analogous of the Poincar\'e--Bendixson theorem for piecewise differential systems. For more details in this extension see for instance \cite{BuzCarEuz2018}.  
\begin{proposition}\cite{CruzNoaesTorre2019}
Let $Z$ be a piecewise differential system as in \eqref{eq:0} with $h(x,y) = y$ such that the origin is a stable monodromic equilibrium. 
Assume $\frac{\partial Y^{+}}{\partial x}(0,0) > 0$. 
Given a real number $b$, let the perturbed system $Z_{b}$ be defined by $Z_b^+=Z^+$ and $Z_b^- = \left(X^-, Y^-+b\right)$. 
Then, for $b$ small enough, the system $Z_{b}$ exhibits a pseudo-Hopf type bifurcation at $b=0$ when $b > 0$. 
See Figure \ref{fi:pseudohopf}. 
\end{proposition}

This is in some sense the canonical form of the bifurcation. Even playing with all  four possible constant terms, it is clear that we will not get more limit cycles. 

As a consequence, in our setting, having found limit cycles of system \eqref{eq:8}, we can always add one more by considering constant terms as in the above proposition. 

\subsection{Expanding the solutions}\label{subsection_aproximation} 
In order to get simpler formulas, we assume from now on that our non-perturbed non-degenerate center $Z = (Z^+, Z^-)$ is written in the canonical form, that is, the linear part of both $Z^{\pm}$ is $(-y, x)$. 
Then the functions $F_i^{\pm}$, $i=0,1,\ldots, m$, of \eqref{eq:drdt} write as 
\begin{equation}\label{F_k}
F^\pm_0(\theta,r)=\dfrac{r^2 \,f_0^\pm(\theta,r)}{1+r\,g_0^\pm(\theta,r)},\quad F^\pm_i(\theta,r)=\dfrac{f_i^\pm(\theta,r)}{\left[1+r\,g_0^\pm(\theta,r)\right]^{l^\pm_i}},
\end{equation}
with $l^\pm_i\in \mathbb{N}^*$, and where $f_0^\pm,g_0^\pm,$ and $f_i^\pm$ are polynomials in $r,\cos\theta$, $\sin\theta$, $a^{\pm}_{\eta},$ and $b^{\pm}_{\eta}$.  

We consider the Taylor expansion in $r$ of $\varphi_i^{\pm}(\theta, r)$:  
\begin{equation}\label{varphi} 
\varphi_i^\pm(\theta,r) = \sum_{j = 1}^N\xi^\pm_{i,j}(\theta)r^j+O\left(r^{N+1}\right), 
\end{equation}
$i = 0,1,\ldots, m$. 
The natural number $N$ is the order of the expansion and functions $\xi^\pm_{i,j}(\theta)$ are suitable ones given by the Taylor series expansion. 
Clearly $\xi^\pm_{0,1} \equiv 1$. 
The idea of the algorithm is to find $\xi^\pm_{i,j}(\theta)$ iteratively as follows. 

We begin with $\varphi_0^{\pm}$. 
By applying the formula of $F_0^{\pm}$, given in \eqref{F_k}, in equation \eqref{eq:drdt}, we get 
\begin{equation*}
\Big(1+\varphi_0^\pm(\theta, r)\,g_0^\pm\left(\theta,\varphi_0^\pm(\theta,  r)\right) \Big){\varphi_0^\pm}'(\theta, r)- \varphi_0^\pm(\theta, r)^2 \,f_0^\pm\left(\theta,\varphi_0^\pm(\theta, r)\right) = 0.
\end{equation*} 
Then by applying the first equation of \eqref{varphi} we recursively get 
\begin{equation*}
\begin{aligned}
{\xi_{0,2}^\pm}'(\theta)& = h^\pm_{0,2}\left(\cos\theta,\sin\theta\right),\\
{\xi_{0,j+1}^\pm}'(\theta)& = h^\pm_{0,j+1}\left(\cos\theta,\sin\theta,\xi_{0,2}^\pm(\theta),\ldots, \xi_{0,j}^\pm(\theta)\right),\\ 
\end{aligned}
\end{equation*}
for $j = 2,\ldots, N-1$, where $h^\pm_{0,j}$'s are suitable \emph{polynomial} functions. 
So, by simple integration of trigonometric polynomials, we obtain the expressions of $\xi_{0,j}^{\pm}$, by using the fact that $\xi_{0,j}^{\pm}(0) = 0$. 
Clearly, if the non-perturbed polynomial center has parameters, they will also appear as variables of the polynomials $h^\pm_{0,j}$'s. 

Now, to obtain the Taylor coefficients $\xi^\pm_{1,j}(\theta)$ of $\varphi_1^{\pm}(\theta, r)$, we consider the differential version of the first formula of \eqref{y1y2y3}, 
\begin{equation*}
{\varphi_1^\pm}'(\theta,r)-\left(F_1^\pm\left(\theta,\varphi_0^\pm(\theta,r)\right) + \frac{\partial F_0^\pm}{\partial r}\left(\theta,\varphi_0^\pm(\theta,r)\right)\varphi_1^\pm(\theta, r) \right) = 0, 
\end{equation*}  
with $\varphi_1^{\pm}(0) = 0$, and the formulas of $F_0^{\pm}$ and $F_1^{\pm}$ given in \eqref{F_k}. 
Then, after applying the second equation of \eqref{varphi}, we recursively get 
\begin{equation*}
\begin{aligned}
{\xi_{1,1}^\pm}'(\theta) & = h^\pm_{1,1}\left(\cos\theta,\sin\theta,a^\pm_{\eta},b^\pm_{\eta}\right),\\ 
{\xi_{1,j+1}^\pm}'(\theta)&=h^\pm_{1,j+1}\left(\cos\theta,\sin\theta,\xi_{1,1}^\pm(\theta),\ldots,\xi_{1,j}^\pm(\theta), a^\pm_{\eta},b^\pm_{\eta}\right), 
\end{aligned}
\end{equation*}
for $j = 1,\ldots, N-1$, where $h^\pm_{1,j}$'s are suitable polynomial functions in $\cos \theta$, $\sin \theta,$ and $a^\pm_{\eta},b^\pm_{\eta}$. 
As above, simple integration provides the expressions of $\xi_{1,k}^{\pm}(\theta)$. 

Acting similarly with the other formulas of \eqref{y1y2y3} and equation \eqref{varphi}, we obtain the coefficients of the Taylor expansions of $\varphi_i^{\pm}$ for $i\geq2$. 
We once more stress that these expressions are obtained in a relatively simple manner, due to the fact that only integration of suitable trigonometric polynomials is necessary. 

Now we define 
\begin{equation}\label{psidefi}
\psi_{i,j} = \xi_{i,j}^+(\pi) - \xi_{i,j}^-(-\pi), 
\end{equation}
for $i = 0,\ldots, m$, $j = 1,\ldots, N$, and the \emph{approximated $i$-difference function} is defined by 
\begin{equation}\label{varphij}
\psi_i(r) = \sum_{j = 1}^N \psi_{i,j} r^j. 
\end{equation}
Since the unperturbed system is a center, it follows that $\psi_{0,j} = 0$ for all $j$. 
Moreover, by the analyticity with respect to the parameters, the $\psi_{i,j}$,for $j = 1, \ldots,$ are homogeneous polynomials of degree $i$ in the parameters, $i = 1, \ldots$ 
By our construction, it is clear that 
$$
\delta_i(r) = \psi_i(r) + \CO(r^{N+1}), 
$$
i.e., the approximated $i$-difference function is the $N$-jet of the $i$-difference function. 

In the smooth case, where the ``$+$ part'' equals the ``$-$ part'' and $\delta_i(\rho) = \varphi_i(2 \pi)$ ($\varphi_i = \varphi_i^+ = \varphi_i^-$), by defining $\psi(r) = \xi(2 \pi)$ (here $\xi = \xi^+ = \xi^-$ as well), it also follows that $\psi_i$ is the $N$-jet of $\delta_i$. 

In case we know that $\delta_i(r) \equiv 0$ for $i=1,\ldots, k-1$ and $\psi_k(r) \neq 0$, then by applying the implicit function theorem it follows that each simple zero of $\psi_k(r)$ close to $r = 0$ provides a simple zero of $\delta_k(r)$, and hence, by a further application of the implicit function theorem, for $\varepsilon> 0$ small enough, each of them provides a simple zero of $r \mapsto \Delta(r, \varepsilon)$. 
Therefore each simple zero of $\psi_k(r)$ will provide a limit cycle of \eqref{eq:8}. 
Considering also the limit cycle coming from the pseudo-Hopf bifurcation, we can summarize the search of limit cycles by looking for zeros of the approximated $k$-difference function in the following result. 

\begin{proposition}\label{hh}
Assume that $\delta_i(r) \equiv 0$, $i=1,\ldots,k-1$ and that $\psi_k(r) \not \equiv 0$. 
If $\psi_k(r)$ has $K$ simple zeros close to $r=0$, then, for all $\varepsilon > 0$ small enough, system \eqref{eq:8} has at least $K + 1$ limit cycles of small amplitude. 
\end{proposition}

This proposition is equivalent to the $k$-order of averaging. 

\begin{remark}
In the smooth case, $K$ simple zeros of $\delta_i(r)$ provides at least $K$ limit cycles of small amplitude. 
\end{remark}

On the other hand, in case we can not assure that $\delta_i(r) \equiv 0$ for $i = 1,\ldots, k-1$, we need to take into account the distinct orders of smallness of $r$ and $\varepsilon$. 
This will be discussed in details in the proofs of Theorems~\ref{theorem:Main2} and \ref{theorem:main}. 

\subsection{Piecewise centers from the Loud family}
The next result analyses the possibilities for system $Z$ of \eqref{eq:0} to have a center by considering $Z^+ = S_i$ and $Z^- = S_j$, with $i \neq j$. 
We use the identification of the straight lines introduced after \eqref{parameter}. 

\begin{lemma}\label{lemma:center}
If $i \neq j$, the system of differential equations $Z$ defined by putting $Z^+ = S_i$ and $Z^- = S_j$ has a center at the origin if and only if $\tau = 0$ or $i,j\in \{1,2\}$. 
\end{lemma}
\begin{proof}
	According to \cite{ChaSab1999}, systems $S_1$, $S_2$, $S_3,$ and $S_4$ have the first integrals $H_1(x,y) = \left(x^2 + y^2\right)/(1 + 2 y)$, $H_2(x,y) = \left(x^2 + y^2\right)/(1 + y)^2$, $H_3(x,y) = \big(9 (x^2 + y^2) - 24 x^2 y + 16 x^4 \big)/(16 y - 3),$ and $H_4(x,y) = \left(9 (x^2 + y^2) + 24 y^3 + 16 y^4\right)/(3 + 8 y)^4$, respectively. 
	
	Given $v = (v_1, v_2) \in \R^2\setminus\{(0,0)\}$, we consider the straight line $L_{v}(\lambda) = \left(v_1 \lambda, v_2 \lambda\right)$, $\lambda\in \R$. 
	For each $i=1,2,3,4$, since $S_i$ has a center at the origin, it follows that for each small enough $\lambda>0$ there exists $\sigma = \sigma_i(\lambda)<0$ such that 
	$$
	H_i(L_v(\lambda)) - H_i(L_v(\sigma)) = 0. 
	$$ 
	Therefore $Z$ will define a center at the origin if and only if $\sigma_i(\lambda) = \sigma_j(\lambda)$ for each $\lambda>0$ small enough. 
	
	If $v_2 = 0$, then $\sigma_i(\lambda) = -\lambda$ for $i=1,2,3,4$. 
	So from now on, we assume $v_2 \neq 0$. 
	
	It is simple to see that $\sigma_1(\lambda) = \sigma_2(\lambda) = -\lambda/(1 + 2 v_2 \lambda)$. 
	
	On the other hand, if $i=3$ or $4$ and $\sigma = \sigma_1(\lambda)$, it is not difficult to conclude that $H_i(L_v(\lambda)) - H_i(L_v(\sigma))$ cannot be identically zero. 
	Hence $\sigma_i(\lambda) \neq \sigma_1(\lambda)$ for $i=3$ or $4$. 
	
	Finally, we consider the numerators of $H_3(L_v(\lambda)) - H_3(L_v(\sigma))$ and $H_4(L_v(\lambda)) - H_4(L_v(\sigma))$ as polynomials in $\sigma$, after canceling the factor $\sigma-\lambda$. 
	If either $v_1\neq 0$ or $v_1 = 0$, the resultant of these polynomials is a non-identically zero polynomial in $\lambda$. 
	So $\sigma_3(\lambda) \neq \sigma_4(\lambda)$ and the lemma is proven. 
\end{proof}

\section{First order non-smooth perturbations of $S_1$, $S_2$, $S_3$, $S_4,$ and $S_1 \& S_2$ for all straight lines}\label{section:S1toS4}
For any given straight line through the origin to be the discontinuity line, we apply our technique in order to study non-smooth perturbations of $Z$ with $Z^{\pm} = S_1$, $Z^{\pm} = S_2$, $Z^{\pm} = S_3$, and $Z^{\pm} = S_4$ in Section~\ref{sub1001}. 
The results on the number of bifurcating limit cycles will prove Theorem~\ref{theorem:order1}. 
Then, in Section~\ref{sub1002}, we do similar study for the discontinuous system $Z$ with $Z^+ = S_1$ and $Z^- = S_2$. 
This will prove Theorem~\ref{theorem:main3}. 

\subsection{First order by using any straight line for each $S_1$--$S_4$}\label{sub1001}
\begin{proof}[Proof of Theorem~\ref{theorem:order1}] 
We begin with $Z = S_1$. 
According to the formulas \eqref{F_k} of Section~\ref{section:averaging}, we have (since $f_0^+ = f_0^-$ and $g_0^{+} = g_0^-$, we simply drop the $\pm$ signals in $f_0$ and $g_0$, also we do not write the parenthesis of the pairs $\eta$) 
\begin{equation}\label{case_S1}
\begin{aligned}
f_0(\theta,r)&=\cos\theta, \quad g_0(\theta,r)=\sin\theta,\\
f^\pm_1(\theta,r)&= \Big( -2 a_{2 0}^{\pm} \cos^3 \theta \sin \theta - 2 a_{11}^{\pm} \cos^2 \theta \sin^2\theta + b_{2 0}^{\pm} \cos^2\theta \big(1 - 2 \sin^2 \theta \big) \\ 
&\phantom{{}=} - \big(b_{1 1}^{\pm} \cos \theta + b_{0 2}^{\pm} \sin \theta \big) \big(1 - 2 \cos^2 \theta \big) \sin \theta - 2 a_{0 2}^{\pm} \cos \theta \sin^3 \theta  \Big)r^3 \\ 
&\phantom{{}=} - \Big(a_{2 0}^{\pm} \cos^3 \theta + \big( b_{1 0}^{\pm} \cos \theta + b_{0 1}^{\pm} \sin \theta \big) \big(1 - 2 \cos^2\theta\big) + \big(a_{1 1}^{\pm} + b_{2 0}^{\pm} \\ 
&\phantom{{}=} + 2 a_{1 0}^{\pm} \big) \cos^2 \theta \sin \theta + \big( 2 a_{0 1}^{\pm} + a_{0 2}^{\pm} + b_{1 1}^{\pm} \big) \cos \theta \sin^2 \theta + b_{0 2}^{\pm} \sin^3 \theta  \Big)r^2 \\ 
&\phantom{{}=} -\Big(\left( a^\pm_{01}+b^\pm_{10}\right) \cos\theta\sin\theta  + a^\pm_{10} \cos^2 \theta + b^\pm_{01} \sin^2 \theta \Big)r,
\end{aligned}
\end{equation}
and $l_1=2$. 
Since we want to consider a generic straight line through the origin making angle $\alpha$ with the $x$-axis, we take the change $\theta\rightarrow \theta+\alpha$ in \eqref{case_S1}, see the beginning of Section~\ref{subsection:averaging}. 
Then, following the steps \eqref{varphi} to \eqref{varphij} of Section~\ref{section:averaging}, we get the recursive expressions below (here we write $\xi_{0,j} = \xi_{0,j}^+ = \xi_{0,j}^-$), for $N = 15$. 
Actually, we do not need such a precision for the first order. 
But since it is required for higher orders, we do the calculations from the very beginning with $N=15$ for using them in Section~\ref{section4} below, for a specific straight line. 
\begin{equation}\label{case_recorrencia:S1}
{\xi_{0,2}}'(\theta)  = \cos(\alpha + \theta), \quad
{\xi_{0,3}}'(\theta)  = - {\xi_{0,2}}'(\theta)\sin(\alpha+\theta) + 2\xi_{0,2}(\theta) \cos(\alpha+\theta), 
\end{equation}
and, similarly, ${\xi_{0,j}}'(\theta)$, $j = 4,\ldots, N$. 
By simple integration we iteratively find the $\xi_{0,j}$'s and, consequently, we get the expression of $\varphi_0(\theta, \rho) = \varphi_0^{\pm}(\theta, \rho)$ according to \eqref{varphi}. 
With this in hands, we continue the algorithm of Section~\ref{section:averaging}, obtaining the equations 
\begin{equation}\label{ffff}
\begin{aligned}
{\xi^{\pm}_{1,1}}'(\theta) & = \left(a^\pm_{10} - b^\pm_{01}\right) \cos^2( \theta+\alpha) + \left(a^\pm_{01} + b^\pm_{10}\right) \cos\left(\theta+\alpha\right)\sin\left(\theta+\alpha\right) + b^\pm_{01}, \\ 
{\xi^\pm_{1,2}}'(\theta)  & = \Big( \left(b^\pm_{01} - a^\pm_{10} \right)\cos^2\left( \theta+\alpha \right) - \left(a^\pm_{01} + b^\pm_{10} \right) \cos\left(\theta+\alpha \right)\sin\left( \theta+\alpha\right) \\ 
&\phantom{{}=} -b^\pm_{01} \Big) \sin \alpha + \left(a^\pm_{20} -2 a^\pm_{01} - a^\pm_{02} - 2 b^\pm_{10} - b^\pm_{11} \right) \cos^3 \left(\theta+\alpha \right) \\ 
&\phantom{{}=} + \left( a^\pm_{10} + a^\pm_{11} - b^\pm_{01} - b^\pm_{02} + b^\pm_{20}\right) \cos^2 \left( \theta+\alpha \right) \sin\left(\theta+\alpha\right) \\ 
&\phantom{{}=} - \left(a^\pm_{01} + b^\pm_{10} \right) \cos\left(\theta+\alpha\right)\sin^2\left(\theta+\alpha \right) + b^\pm_{02} \sin \left(\theta+\alpha \right) \\ 
&\phantom{{}=} + \left(2 a^\pm_{01} + a^\pm_{02} + b^\pm_{10} + b^\pm_{1,} \right) \cos \left(\theta+\alpha \right) +2\cos\left(\theta+\alpha \right)\xi^\pm_{1,1}(\theta), 
\end{aligned}
\end{equation}
and ${\xi^{\pm}_{1,j}}'(\theta)$ for $j = 3,\ldots, N$. 
These equations can also be solved iteratively by simple integration. 
But before doing this we write them as polynomials in $\cos \alpha$ and $\sin \alpha$ and use the rational parametrization of $(\cos \alpha, \sin \alpha)$ given by \eqref{parameter}, and so in particular we introduce the notation of the theorem. 
Then we integrate the obtained equations and, by considering \eqref{psidefi}, we get 
\begin{equation*}
\begin{aligned}
\psi_{1,1}& =  \pi\left(a^+_{10}+b^+_{01}+a^-_{10}+b^-_{01} \right)/2, \\
\psi_{1,2} & = \Big( 16\tau^3(a^-_{02} -a^+_{02} +a^+_{01}-a^-_{01} -b^+_{11}+b^-_{11}) +4 \tau\big(3\tau^4+2\tau^2+3\big) (a^-_{20} - a^+_{20}) \\
&\phantom{={}}  -\big(2 (\tau^6 - 3 \tau^4 +3 \tau^2 -1) +9\pi (\tau^5 + 2 \tau^3 + \tau) \big) a^-_{10} +\big(2 (\tau^6 -3\tau^4 +3\tau^2 -1)  \\ 
&\phantom{={}}  - 9\pi (\tau^5 + 2 \tau^3 +\tau) \big)a^+_{10} -\big( 2 (5\tau^6 +9 \tau^{4} -9\tau^2 -5) +9\pi(\tau^5 + 2 \tau^3  + \tau)  \big) b^-_{01}  \\
&\phantom{={}} + \big( 2(5\tau^6 +9\tau^4 - 9\tau^2 -5) - 9\pi (\tau^5 +2 \tau^3 + \tau) \big) b^+_{01} +4\tau\big(3\tau^4+10\tau^2+3\big) (b^+_{10}  \\
&\phantom{={}} - b^-_{10}) +2\big(\tau^6-3\tau^{4}+3\tau^2-1\big) (a^-_{11} - a^+_{11})  +4\big(\tau^6 +3\tau^4-3\tau^2-1 \big) (b^-_{02} - b^+_{02}) \\
&\phantom{={}} + 2\big(\tau^6-3\tau^4+3\tau^2 -1\big) (b^-_{20} - b^+_{20}) \Big)/\left(3(\tau^2+1)^3\right) 
\end{aligned}
\end{equation*} 
and $\psi_{1,j}$, $j = 3,\ldots, N$, with $\tau\in[-1,1)$. 

We write $\psi_1^\tau(r) = \psi_1(r) = \sum_{j = 1}^{N}\psi_{1,j} r^j$ to express the dependency on $\tau$. 
Then noting that the $\psi_{i,j}$'s, are linear in the parameters $a^\pm_{\eta}$ and $b^\pm_{\eta}$, $|\eta|\leq 2$, we introduce new variables $\alpha_j$ by iteratively solving in $a^\pm_{\eta}$ or $b^\pm_{\eta}$ the equations $\psi_{1,1} = \alpha_1$, $\psi_{1,2} = \alpha_2$ and so on, until we get a that $\psi_{1,j}$ depends only on $\alpha_l$, $l = 1,2, \ldots, j-1$. 
Then we write $\widetilde{\alpha}_j = \psi_{1,j}$ and continue this process for $j+1$, $j+2,\ldots,N$. 
That is, we rewrite $\psi_1^\tau$ for different $\tau$'s as 
\begin{equation*}
\begin{aligned}
\psi^\tau_1(r)&=\sum_{j = 1}^{4}\alpha_j r^j+\widetilde{\alpha}_5 r^5 +\alpha_6 r^6+\widetilde{\alpha}_7 r^7+\alpha_8 r^8 +\widetilde{\alpha}_9 r^9 + \alpha_{10} r^{10} \\ 
& \phantom{={}} +\sum_{j = 11}^{N} \widetilde{\alpha}_j r^{j},\  \tau\in (-1,1)\setminus\{0\} , \\
\psi^\tau_1(r) & = \sum_{j = 1}^{4}\alpha_j r^j + \widetilde{\alpha}_5 r^5 + \alpha_6 r^6 + \sum_{j = 7}^{N}\widetilde{\alpha}_j r^j,\  \tau= -1, 0,
\end{aligned}
\end{equation*}
where $\alpha_j$ (depending on the parameters $a^\pm_{\eta}$ and $b^\pm_{\eta}$) can be made any real number, and $\widetilde{\alpha}_j$ depends on $\alpha_{l}$ for $l < j$ in such a way that 
$$
\alpha_l = 0,\ l = 1,\dots,j - 1 \implies \widetilde{\alpha}_j = 0, 
$$
for $j = 1,2,\ldots, N$. It is then easy to find parameters $\alpha_j$'s such that $\psi_1^\tau$ have $6$ (respectively $4$) simple positive zeros if $\tau\in (-1,1)\setminus\{0\}$ (respectively $\tau= -1,0$). Actually, in order to obtain these zeros, we even do not need the full freedom of $\alpha_{10}$ (respectively $\alpha_6$), it is only necessary that if $\alpha_i = 0$, $i = 1, \ldots, 9$, then $\alpha_{10} \neq 0$ (respectively $\alpha_i = 0$, $i = 1, \ldots, 5$, then $\alpha_6 \neq 0$). Then we apply Proposition~\ref{hh} to conclude that for $\tau \in(-1,1)\setminus \{0\}$ we get $7$ limit cycles, while for $\tau=-1$ or $\tau=1$ we get $5$ limit cycles. This proves the theorem for $S_1$. 

For $S_2$ the calculations are completely analogous, so we do not detail them here. 
We just write down the approximated averaged functions for the different $\tau$'$s$ as above: 
\begin{equation*}
\begin{aligned}
\psi^\tau_1(r)& = \sum_{j = 1}^{6}\alpha_j r^j+\widetilde{\alpha}_7 r^7+\alpha_8 r^8 +\widetilde{\alpha}_9 r^9 + \alpha_{10} r^{10} + \sum_{j = 11}^{N} \widetilde{\alpha}_j r^{j},\  \tau \in(-1,1)\setminus\{0\}, \\
\psi^\tau_1(r) & = \sum_{j = 1}^{6}\alpha_j r^j+\sum_{j = 7}^{N}\widetilde{\alpha}_j r^j,\  \tau=0,1,
\end{aligned}
\end{equation*}
with $\alpha_j$'s and $\widetilde{\alpha}_j$'s as in the $S_1$ case. 
Therefore, as above, we can find at least $7$ (respectively $5$) zeros of $\psi_1^\tau$ and hence Proposition~\ref{hh} gives at least $8$ (respectively $6$) limit cycles if $\tau \in (-1, 1)\setminus\{0\}$ (respectively $\tau=0,1$), proving the theorem for $S_2$. 

For $S_3$ and $S_4$, the calculations are also analogous. 
The results are as stated in the theorem. 
The jumps of independence here are a bit different, as they not only change when $\tau=0$ or $1$, but also when we have $\tau$ as real solutions of some suitable polynomials. 
Anyway, the number of independent coefficients does not change for these special $\tau$'s. 
\end{proof}

\subsection{First order for $Z$ with $Z^+ = S_1$ and $Z^- = S_2$, for any straight line}\label{sub1002}
\begin{proof}[Proof of Theorem~\ref{theorem:main3}] Lemma~\ref{lemma:center} guarantees that the origin is a center. Since $Z^+ = S_1$, the expressions of $f^+_0, g^+_0,$ and $f^+_1$ are given in \eqref{case_S1}. The expressions for $Z^- = S_2$ are 
\begin{equation}\label{case_S2_fg}
\begin{aligned}
f^-_0(\theta,r)&=\cos\theta,\quad g^-_0(\theta,r)=0, \\
f^-_1(\theta,r)&=\Big(-b_{2 0} \cos^4 \theta + \big(a_{2 0}^--b_{1 1}^- \big) \cos^3 \theta \sin \theta +\big(a_{1 1}^--b_{0 2}^- \big) \cos^2 \theta \sin^2 \theta \\ 
& \phantom{={}} + a_{0 2}^- \cos \theta \sin^3 \theta \Big)r^3 +\Big( \big( a_{2 0}^--b_{1 0}^- \big)\cos^3 \theta + \big( a_{1 0}^-+a_{1 1}^--b_{0 1}^- \\ 
& \phantom{={}} +b_{2 0}^-  \big) \cos^2 \theta \sin \theta + \big( a_{0 1}^- + a_{0 2}^-+b_{1 1}^-\big) \cos \theta \sin^2 \theta + b_{0 2} \sin^3 \theta  \Big)r^2 \\ 
& \phantom{={}} + \Big( a^-_{10} \cos^2\theta + \left(a^-_{01} +b^-_{10}\right)\cos\theta \sin\theta + b^-_{01} \sin^2 \theta\Big)r, 
\end{aligned}
\end{equation}
and $l_1^-=2.$	We take the change $\theta\rightarrow \theta+\alpha$ in \eqref{case_S1} and \eqref{case_S2_fg}. 
Then following the steps \eqref{varphi} to \eqref{varphij} of Section~\ref{section:averaging} with $N = 15$, we get (for $Z^+ = S_1$) the $\xi_{0,j}^+$'s in \eqref{case_recorrencia:S1}. 
And for $Z^- = S_2$ it follows that  
\begin{equation*}
\begin{aligned}	
{\xi_{0,2}^-}'(\theta)  =& \cos(\alpha + \theta), \quad {\xi_{0,3}^-}'(\theta)  = 2\xi^-_{0,2}(\theta) \cos(\alpha+\theta),\\ 
{\xi_{0,4}^-}'(\theta)  =& \left(2\xi^-_{0,3}(\theta)+{\xi^-_{0,2}}^2(\theta)\right) \cos(\alpha+\theta),
\end{aligned}
\end{equation*}
and similarly ${\xi_{0,j}^-}'$ for $j = 5, \ldots, N$. 
By simple integration we find the $\psi^\pm_{0,i}$'s and, consequently, we get the expression of $\varphi_0^{\pm}(\theta, r)$. 
With this in hands, we continue the algorithm of Section~\ref{section:averaging}, obtaining the equations \eqref{ffff} for $Z^+ = S^1$, and, for $Z^- = S_2$ we get 
\begin{equation*}\label{gggg}
\begin{aligned}
{\xi^{-}_{1,1}}'(\theta) & = (a^-_{01} + b^-_{10})\cos(\theta+\alpha)\sin( \theta+\alpha )+ (a^-_{10} - b^-_{01} ) \cos^{2}(\theta+\alpha) +  b^-_{01},\\
{\xi^-_{1,2}}'(\theta)  & = \Big( (b^-_{01} - a^-_{10}) \cos^{2} ( \theta +\alpha)  - (a^-_{02} + b^-_{10}) \cos( \theta +\alpha)\sin(\theta +\alpha) \\ 
& \phantom{={}} - b^-_{01} \Big) \sin(\alpha) + (a^-_{20} - a^-_{01} - a^-_{02} - b^-_{10} - b^-_{11})\cos^{3} ( \theta +\alpha) \\ 
& \phantom{={}} + (2 a^-_{10} +a^-_{11} -b^-_{02} + b^-_{20} -2 b^-_{01}) \cos^{2} ( \theta +\alpha  )\sin(\theta +\alpha) \\
& \phantom{={}} + (a^-_{01} +b^-_{10} ) \cos(\theta +\alpha )\sin^{2}( \theta +\alpha) + (b^-_{01} + b^-_{02}) \sin( \theta +\alpha) \\
& \phantom{={}} + (a^-_{01} + a^-_{02} + b^-_{11}) \cos(\theta +\alpha)  + 2\cos(\theta +\alpha)\xi^{-}_{1,1} ( \theta), 
\end{aligned}
\end{equation*}
and ${\xi^-}_{1,j}'(\theta)$ for $j = 3,\ldots, N$. 
As above, these can be iteratively solved by simple integration. 
But before doing this, as above, we introduce the parametrization of $(\cos \alpha, \sin \alpha)$ given by \eqref{parameter}. 
Then we integrate and, by considering \eqref{psidefi}, we get 
\begin{equation}\label{ordemum_S1S2}
\begin{aligned}
\psi_{1,1} & = \pi\left(a^+_{10} + b^+_{01} + a^-_{10} + b^-_{01} \right)/2,  \\
\psi_{1,2}  & =  \Big(16 \tau^{3} (a^+_{02} - a^-_{02} - b^-_{11} + b^+_{11} - a^+_{01}) + 9\pi \left(\tau^{5} + 9\tau^{3} + \tau \right) a^-_{10} + \big(9 \pi (\tau^5  \\
& \phantom{={}} + 2 \tau^{3} + \tau) - 2 (\tau^{6} - 3 \tau^{4} + 3 \tau^{2} - 1)\big) a^+_{10} + 4 \tau\big(3 \tau^{4} + 2 \tau^{2} + 3 \big) (a^+_{20} - a^-_{20}) \\
& \phantom{={}} +2\left(\tau^{6} - 3 \tau^{4} + 3 \tau^{2} - 1\right) (a^+_{11} - a^-_{11}) +3\big(3 \pi (\tau^{5} + 2 \tau^{3} + \tau) + 2 (\tau^{6} + \tau^{4} \\
& \phantom{={}} - \tau^{2} - 1)\big) b^-_{01} +\big(9 \pi (\tau^{5} + 2 \tau^{3} + \tau) - 2 (5 \tau^{6} + 9 \tau^{4} - 9 \tau^{2} - 5) \big) b^+_{01} - 4\big(\tau^{6}  \\
& \phantom{={}}  + 3 \tau^{4} - 4\tau^{2} +4\big) b^-_{02} + 4\big( \tau^{6} + 3 \tau^{4} -3 \tau^{2} - 4 \big) b^+_{02} - 4 \tau\big(3\tau^{4} + 10 \tau^{2} + 3\big) b^+_{10} \\
& \phantom{={}}  + 12 \big(\tau^{5} + 2 \tau^{3} + \tau \big) b^-_{10} + 2\big(\tau^{6} - 3 \tau^{4} +3 \tau^{2} - 1\big) (b^+_{20} - b^-_{20}) \Big)/\big(3(\tau^2 + 1)^3\big) 
\end{aligned}
\end{equation} 
and similarly $\psi_{1,j}$, $j = 3,\ldots, N$, with $\tau\in[-1,1)$. 
We write $\psi_1^\tau(r) = \psi_1(r) = \sum_{j = 1}^{14}\psi_{1,j} r^j$ to express the dependency on $\tau$. 
Then, as Section~\ref{sub1001}, we iteratively solve $\psi_{1,j} = \alpha_j$ (when it is possible) and rewrite $\psi_1^\tau(r)$ as 
\begin{equation*}
\begin{aligned}
\psi^\tau_1(r)&=\sum_{j = 1}^{8}\alpha_j r^k + \widetilde{\alpha}_9 r^9 + \alpha_9 r^{10} + \widetilde{\alpha}_{11} r^{11} + \alpha_{10} r^{12} +\sum_{j = 13}^{N} \widetilde{\alpha}_j r^{j},\  \tau\in(-1,1)\setminus \{0\}, \\
\psi^\tau_1(r) & = \sum_{j = 1}^{6}\alpha_j r^j + \sum_{j= 7}^{N}\widetilde{\alpha}_j r^j,\  \tau = 0,  \hspace{1cm}
\psi^\tau_1(r)  = \sum_{j = 1}^{8}\alpha_j r^j + \sum_{j = 9}^{N}\widetilde{\alpha}_j r^k,\  \tau=-1, 
\end{aligned}
\end{equation*}
where as above $\alpha_j$ (depending on the parameters $a^\pm_{\eta}$ and $b^\pm_{\eta}$) can be made any real number, and $\widetilde{\alpha}_j$ depends on $\alpha_l$ for $l<j$ in such a way that 
$\alpha_l = 0,\ l=1,\dots,j-1$ implying that $\widetilde{\alpha}_j = 0$, for $j=1,\ldots, N$. 
Hence, as above, after applying Proposition~\ref{hh}, it follows that for any $\tau\in(-1,1)\setminus\{0\}$, we find $10$ limit cycles. 
Analogously, for $\tau = 0$ and for $\tau= -1$, there are $6$ and $8$ limit cycles, respectively. 
\end{proof}

\section{Higher order non-smooth perturbations of $S_1$, $S_2$, $S_3$, $S_4$, and $S_1 \& S_2$ for a specific straight line}\label{section4}
\subsection{Perturbations of $S_1$, $S_2$, $S_3$, and $S_4$}	
\begin{proof}[Proof of Theorem~\ref{theorem:Main2}] 
In this proof, we will provide the calculations for $S_4$. Systems $S_1$, $S_2$, and $S_3$ have analogous reasoning. Indeed, for $S_1$ and $S_2$ almost the same algorithm applies. For $S_3$, there is a slight difference, that we anyway point out after the end of the proof for $S_4$. We assume the notations of Section~\ref{section:averaging}. We already know by Theorem~\ref{theorem:order1} that for first order, the averaging theory provides at least $9$ limit cycles. 
Here we keep $N = 15$. In the situation of $S_4$ given in \eqref{eq:10}, we have, according to \eqref{F_k}: $f^\pm_0(\theta,r) = 12\cos\theta\left(3\cos^2\theta+1\right)$, $g_0(\theta,r)=4\sin^3\theta -8/3\sin\theta$, 
\begin{equation*}
\begin{aligned}
f^{\pm}_1(\theta,r) & = 12\Big(-2 a^\pm_{2 0}\cos^3\theta \sin\theta  - 2 a^\pm_{1 1}\cos^2\theta \sin^2\theta + \big( b_{2 0}^{\pm} \cos^2\theta + b_{1 1}^{\pm} \cos \theta \sin \theta  \\ 
& \phantom{={}} + b^\pm_{02} \sin^2 \theta  \big) \big(5 \cos^2 \theta - 1\big)  -2a^\pm_{02} \cos\theta \sin^3\theta \Big)r^3 + 3\Big( -3 a^\pm_{2 0} \cos^3\theta  \\
& \phantom{={}} + 4 \big( b^\pm_{10} \cos\theta + b^\pm_{01} \sin \theta \big) \big(5\cos^2\theta - 1\big) - \big(8 a^\pm_{10} + 3 a^\pm_{11} + 3 b^\pm_{20}  \big) \cos^2\theta \sin\theta  \\ 
& \phantom{={}} - \big( 8 a^\pm_{01} + 3 a^\pm_{02} + 3 b^\pm_{11} \big) \cos \theta \sin^2\theta  -3 b^\pm_{02} \sin^3\theta \Big)r^2  - 9\Big( a^\pm_{10} \cos^2\theta  \\
& \phantom{={}} + \big(a^\pm_{01} + b^\pm_{10} \big) \cos\theta \sin\theta + b^\pm_{01} \sin^2\theta  \Big)r 
\end{aligned}
\end{equation*} 
and $f_2^{\pm}(\theta,r)= \sum_{j=1}^{4}\Gamma^\pm _{2,j}\left(\cos\theta,\sin\theta,a^\pm_{\eta},b^\pm_{\eta}\right)r^j,$ where $\left |  \eta\right |\leq2 $, and $\Gamma^\pm_{2,j},$ are homogeneous polynomials of degree two in the perturbative parameters. Here $l_k=k+1$, with $k=0,1,2$. Due to the size of the expressions of $f_2^{\pm}$, we will not present them explicitly.

Here we are assuming that the line of discontinuity is given by $\tau=1/2$, i.e., $\alpha$ is such that $(\cos \alpha, \sin \alpha) = (3/5, 4/5)$. 
By considering the change of variables given by $\theta \mapsto \theta + \alpha$, simple trigonometric relations prove that this is the same as changing $\cos \theta$, $\sin \theta$ by 
\begin{equation}\label{trans_S2}
\frac{3}{5}\cos\theta -\frac{4}{5}\sin\theta, \quad    \frac{4}{5}\cos\theta + \frac{3}{5}\sin\theta, 
\end{equation}
respectively. 
By following the algorithm of Section~\ref{subsection_aproximation}, we get $\psi_1(r) = \sum_{j = 1}^{N}\psi_{1,j} r^j$, with the first three coefficients given by 
\begin{equation*}
\begin{aligned}
\psi_{1,1} & = \pi\left(a^+_{10}+b^+_{01}+a^-_{10}+b^-_{01} \right)/2,  \\ 
\psi_{1,2} & =\Big(1024 (a^-_{01} - a^+_{01}) + 9600 (a^-_{02} - a^+_{02}) + 4050 (a^+_{11} - a^-_{11}) + 35400 (a^-_{20} - a^+_{20}) \Big.\\
& \phantom{={}} + 200576 (b^+_{10} - b^-_{10}) -\left(151200 \pi +158832\right) b^+_{01} -\left(151200 \pi -158832\right) b^-_{01} \\
& \phantom{={}}  -\left(151200 \pi -10368\right) a^-_{10} -\left(151200 \pi +10368\right) a^+_{10} + 9600 (b^-_{11} - b^+_{11})  \\
& \phantom{={}}  + 29700 (b^+_{02} - b^-_{02}) + 4050 (b^+_{20} - b^-_{20}) \Big)/28125, \\
\psi_{1,3}& = \Big(225 \big( (15625\pi - 64512)  a^+_{11} + (15625\pi + 64512) a^-_{11} \big) + 34406400 (a^+_{02} - a^-_{02} \\
& \phantom{={}} + b^+_{11} - b^-_{11}) + 126873600 (a^+_{20} - a^-_{20})  + 8 \big( (44939675\pi - 71156736) b^-_{01}   \\
& \phantom{={}} + (44939675\pi + 71156736) b^+_{01} \big) +3670016 (a^+_{01} - a^-_{01}) + 2048\big( (183175 \pi  \\
& \phantom{={}} - 18144) a^-_{10} + (183175 \pi + 18144) a^+_{10}  \big) - 300 \big( (15625\pi - 354816)b^-_{02}  \\
& \phantom{={}} + (15625\pi + 354816 )b^+_{02} \big)  + 718864384 (b^-_{10} - b^+_{10}) - 150 \big( (78125\pi \\
& \phantom{={}} - 96768) b^-_{20} + (78125 \pi + 96768) b^+_{20} \big) \Big)/7031250. 
\end{aligned}
\end{equation*}
(We explicit only three coefficients because of the big size of the others.) 
Now we calculate $\psi_2(r) = \sum_{j = 1}^{N}\psi_{2,j} r^j$ following the algorithm of Section~\ref{subsection_aproximation}. 
The coefficients $\psi_{2,j}$ are 
\begin{equation*}\label{coef_o2}
\begin{aligned}
\psi_{2,1} & = \pi^2 \left({a^+_{10}}^{2} - {a^-_{10}}^2 + {b^+_{01}}^{2} - {b^-_{01}}^{2}\right) + 2 \pi \Big( b^+_{01} b^+_{10} - \left(a^-_{10} + b^-_{01}\right) a^-_{01}  \\
& \phantom{={}} + \left(\pi b^+_{01} + b^+_{10}\right) a^+_{10} - \left( \pi b^-_{01}- b^-_{10}\right) a^-_{10} - \left( a^+_{10}+ b^+_{01}  \right) a^+_{01} + b^-_{01} b^-_{1 0} \Big) 
\end{aligned}
\end{equation*}
and $\psi_{2,j} =\Omega_{2,j}\left(a^\pm_{\eta},b^\pm_{\eta}\right)$, $j = 2, \ldots, N$, where $\Omega_{2,j}$'s, $j = 2,\ldots, N$ are suitable homogeneous polynomials of degree $2$ in $a^\pm_{\eta}, b^\pm_{\eta}$, $|\eta|\leq 2$. 

Here we cannot apply Proposition~\ref{hh}, because even annihilating $\psi_1$, we do not know whether $\delta_1\equiv 0$. 
In order to proceed, we consider the complete expansions of the difference function $\Delta(r, \varepsilon)$ in $\varepsilon$ and $r$. 
This is an analytic function in $(r,\varepsilon)$, hence we can write: 
$$
\Delta(r, \varepsilon) = \sum_{i=1}^{\infty} \varepsilon^i \sum_{j = 1}^{\infty} \psi_{i,j} r^j = \sum_{j = 1}^{\infty} \left( \sum_{i = 1}^{\infty} \varepsilon^i \psi_{i,j}\right) r^j. 
$$
Each ``coefficient'' $\psi_{i,j}$ is a homogeneous polynomial of degree $i$ in $a^{\pm}_{\eta}$, $b^{\pm}_{\eta}$, $i, j = 1, 2, \ldots$ 
Now we ``eliminate'' $\varepsilon$ by redefining the perturbative coefficients as 
\begin{equation}\label{elimina}
\widetilde{a}^\pm_{\eta} = \varepsilon a^\pm_{\eta}, \quad \widetilde{b}^\pm_{\eta} = \varepsilon b^\pm_{\eta}, 
\end{equation} 
for $|\eta|\leq 2$. 
That is, the perturbation now is such that all the coefficients are ``small''. 
It follows that the difference function now does not depend on $\varepsilon$ anymore: 
$$
\Delta(r) = \sum_{j = 1}^{\infty} \left(\sum_{i = 1}^{\infty} \widetilde{\psi}_{i,j} \right) r^j, 
$$ 
where $\widetilde{\psi}_{i,j}$ is a homogeneous polynomial of degree $i$ in the variables $\widetilde{a}^{\pm}_{\eta}$, $\widetilde{b}^{\pm}_{\eta}$, $i,j = 1, 2, \ldots$ 

We remark that in our algorithm up to now we have explicitly calculated the coefficients $\widetilde{\psi}_{i,j}$ for $i = 1,2$ and $j = 1, \ldots, N$. 
Below we will argue with them in order to find suitable perturbative coefficients guaranteeing that $\Delta(r)$ has at least a certain number $K$ of zeros close to $r=0$, and so system \eqref{eq:8} will have at least $K+1$ limit cycles (by adding the extra limit cycle coming from the pseudo-Hopf bifurcation according to Section~  \ref{pseudo}). 
We define 
$$
\widetilde{\Psi}_j = \sum_{i=1}^{\infty} \widetilde{\psi}_{i,  j},  
$$
$j = 1, \ldots, N$. 
We first analyze the linear part of each $\Psi_j$, namely, $\widetilde{\psi}_{1,j}$. 
Similarly as we have done in the preceding section, we begin rewriting the $\widetilde{\psi}_{1,j} = \alpha_j$ solving the appearing linear equations for $j = 1$, then $j = 2$ and so on until we find a coefficient $\widetilde{\psi}_{1,j}$ depending only on the $\alpha_l$, $l = 1,\ldots, j-1$. 
In our case here we get $\widetilde{\psi}_{1,j} = \alpha_j$ for $j = 1,\ldots, 6$, solving in $\widetilde{a}^-_{1 0}$, $\widetilde{b}^-_{0 1}$, $\widetilde{a}^-_{1 1}$, $\widetilde{a}^-_{0 2}$, $\widetilde{b}^-_{2 0}$, $\widetilde{a}^-_{2 0}$, respectively. 
Then $\widetilde{\psi}_{1,7}$ depends only on $\alpha_{l}$, $l = 1,\ldots, 6$. 
We can thus solve $\widetilde{\psi}_{1,8} = \alpha_7$ in $\widetilde{a}^+_{1 1}$. 
Then it turns out that $\widetilde{\psi}_{1,9}$ depends only on $\alpha_{l}$, $l = 1,\ldots, 7$. 
Further we solve $\widetilde{\psi}_{1,10} = \alpha_8$ in $\widetilde{b}^+_{2 0}$. 
Then $\widetilde{\psi}_{1,11}$ depends only on $\alpha_{l}$, $l = 1,\ldots, 8$. 
Finally, we solve $\widetilde{\psi}_{1,12} = \alpha_9$ in $\widetilde{b}^+_{0 1}$, getting that $\widetilde{\psi}_{1,13}$, $\widetilde{\psi}_{1,14}$ and $\widetilde{\psi}_{1,15}$ all depend only on $\alpha_{l}$, $l = 1, \ldots, 9$. 
Precisely, we succeeded to rewrite the $\widetilde{\psi}_{j,1}$'s as 
\begin{equation*}\label{linear_o1}
\begin{aligned}
\widetilde{\psi}_{1,l} & =\alpha_l,\ l = 1,\ldots, 6,\quad \widetilde{\psi}_{1,8} = \alpha_7, \quad \widetilde{\psi}_{1,10} = \alpha_{8} \quad \widetilde{\psi}_{1,12} = \alpha_{9},\\
\widetilde{\psi}_{1,7} & =-\tfrac{13061776996188618752 }{2780914306640625} \alpha_1 -\tfrac{54153241671606272 }{7415771484375}\alpha_2 -\tfrac{1319866958176 }{263671875}\alpha_3\\
& \phantom{={}}  -\tfrac{92382032896 }{52734375}\alpha_4-\tfrac{44894744 }{140625}\alpha_5-\tfrac{3584 }{125}\alpha_6,\\
\widetilde{\psi}_{1,9} &=\tfrac{35427806878368205783957504 }{14483928680419921875}\alpha_1 +\tfrac{1313963570140369073668096 }{347614288330078125}\alpha_2 \\
& \phantom{={}} +\tfrac{7121695391403890212096 }{2780914306640625}\alpha_3 +\tfrac{1296197793646163968 }{1483154296875}\alpha_4\\
& \phantom{={}}+\tfrac{2960772455199952 }{19775390625}\alpha_5+\tfrac{23193321472 }{2109375}\alpha_6 -\tfrac{896 }{25}\alpha_7,\\
\widetilde{\psi}_{1,11} & = -\tfrac{30988402760095390237763808174014464 }{18331222236156463623046875}\alpha_1 -\tfrac{382824149846431286161826196488192}{146649777889251708984375}\alpha_2 \\
&\phantom{={}} -\tfrac{76744035193724040416612082688}{43451786041259765625}\alpha_3 -\tfrac{69650183483775462018056192 }{115871429443359375}\alpha_4 \\
&\phantom{={}} -\tfrac{284726400765915795149312 }{2780914306640625}\alpha_5 -\tfrac{18339444472611340288 }{2471923828125}\alpha_6  \\ 
&\phantom{={}} + \tfrac{354937470976 }{17578125}\alpha_7 -\tfrac{5376 }{125}\alpha_8,
\end{aligned}
\end{equation*}
and $\widetilde{\psi}_{1,j}$, $j=13,14,15$, also depending only on $\alpha_l$, $l = 1\ldots, 9$. 
We remark that at least $8$ simple zeros of $\Delta(r)$ are already guaranteed with the calculations until here. 
Indeed, for instance, we can reintroduce the $\varepsilon$ and apply Proposition~\ref{hh}. 
(In particular, this agrees with the $9$ limit cycles of Theorem~\ref{theorem:order1} for $\tau = 1/2$.) 

We also remark that if we do not have the ``jumps of independence'' in $\widetilde{\psi}_{1,7}$, $\widetilde{\psi}_{1,9}$ and $\widetilde{\psi}_{1,11}$, that is, if we could write $\widetilde{\psi}_{1,j} = \alpha_j$ $j=1,\ldots, 12$ by means of a linear change of variables, then we could set all the other perturbative coefficients to zero and, by using the implicit function theorem successively for $j = 1, \ldots, 12$, we can solve $\Psi_j = \beta_j$ for given $\beta_j$ in a certain small interval around $0$, getting analytic functions $\alpha_j$ in the variables $\beta_1, \beta_2, \ldots, \beta_j, \alpha_{j+1}, \ldots, \alpha_{12}$. 
Then, after setting all the other parameters to zero, the functions $\Psi_j$, $j > 12$, would be analytic ones in the variables $\beta_1, \ldots, \beta_{12}$ with $\Psi_j(0) = 0$, and we would be able to write 
$$
\Delta(r) = \Delta_{\beta}(r) = r \beta_1\left(1 + \CO(r^{12}) \right) + r^2 \beta_2 \left( 1 + \CO(r^{11}) \right) + \cdots + r^{12} \beta_{12} \left(1 + \CO(r)\right), 
$$ 
for free small enough $\beta_1,\ldots,\beta_{12}$. 
Here we add the subscript $\beta$ to $\Delta$ to emphasize that we are restricting the function $\Delta(r)$, which (also) depends on a priori of all the perturbative parameters to the subspace of only $12$ parameters $\beta = (\beta_1, \ldots, \beta_{12})$, all the others set to zero. 
It would be then simple to obtain at least $11$ small zeros of $\Delta(r)$, so we would get at least $12$ limit cycles of system \eqref{eq:8} (according to Section~\ref{pseudo}). 
We observe that here we do not need the complete independence (on $\beta_1, \ldots, \beta_{11}$) of $\beta_{12}$. 
When it depends continuously on $\beta_1, \ldots, \beta_{11}$, it is enough that it is non-zero when $\beta_1 = \cdots = \beta_{11} = 0$. 
We further note that since $\Delta_{\beta}(r)$ is restricting $\Delta(r)$ to a subspace with some of the free parameters set to zero, the complete function $\Delta(r)$ could have more zeros than $\Delta_{\beta}(r)$. 
Therefore, we could not talk about upper bounds of the number of limit cycles with this technique. 

But we do have the ``jumps'', and we can not apply the implicit function theorem directly. 

Anyway, the use of the implicit function theorem in the preceding argument is equivalent to considering the analytic map $\left(\Psi_1, \Psi_2, \ldots, \Psi_{11}\right): \R^{11} \to \R^{11}$ and to prove that it is a local diffeomorphism at $(\alpha_1, \ldots, \alpha_{11}) = 0$ (all the other variables set to zero), so that $\Psi_{12}$ is not zero when $(\alpha_1, \ldots, \alpha_{11}) = 0$. 

We will be able to do this by considering also the parts of order $2$, namely, $\widetilde{\psi}_{2,j}$ of the $\Psi_j$'s with ``dependent'' linear parts and by considering a special blow up in the parameters. 

Indeed, we act in two steps. 
First we make an appropriate linear change of variables in order to eliminate the linear terms from $\widetilde{\Psi}_7$, $\widetilde{\Psi}_{9},$ and $\widetilde{\Psi}_{11}$. 
Then in particular the new $6$ functions (that we anyway keep the notation) $\widetilde{\Psi}_7,\ldots,\widetilde{\Psi}_{12}$ have no linear terms in the variables $\alpha_1,\ldots,\alpha_6$, and so we do not loss generality in assuming $\alpha_j = 0$ for $j = 1,\ldots,6$ (this will simplify the computations). 
Then, in order to work with only $6$ variables, we consider the simplification 
$$
\widetilde{a}_{0 1}^- = \widetilde{a}_{0 1}^+ = \widetilde{a}_{10}^+ = \widetilde{a}_{2 0}^+ = \widetilde{b}_{0 2}^- = \widetilde{b}_{1 0}^- = \widetilde{b}_{1 1}^- = \widetilde{b}_{0 2}^+ = 0. 
$$
Therefore we end up with $6$ analytic functions in the variables $\alpha_7$, $\alpha_8$, $\alpha_9$, $\widetilde{a}_{0 2}^+$, $\widetilde{b}_{1 0}^+,$ and $\widetilde{b}_{1 1}^+$: 
\begin{equation}\label{Psi_1}
\begin{aligned}
\widetilde{\Psi}_7 & = g_{7}(\alpha_7,\alpha_8,\alpha_9,\widetilde{a}_{0 2}^+, \widetilde{b}_{1 0}^+, \widetilde{b}_{1 1}^+), & 
\widetilde{\Psi}_8 & = \alpha_7 + g_{8}(\alpha_7,\alpha_8,\alpha_9,\widetilde{a}_{0 2}^+, \widetilde{b}_{1 0}^+, \widetilde{b}_{1 1}^+), \\
\widetilde{\Psi}_9 & = g_{9}(\alpha_7,\alpha_8,\alpha_9,\widetilde{a}_{0 2}^+, \widetilde{b}_{1 0}^+, \widetilde{b}_{1 1}^+), & 
\widetilde{\Psi}_{10} & = \alpha_{8} + g_{10}(\alpha_7,\alpha_8,\alpha_9,\widetilde{a}_{0 2}^+, \widetilde{b}_{1 0}^+, \widetilde{b}_{1 1}^+), \\
\widetilde{\Psi}_{11} & = g_{11}(\alpha_7,\alpha_8,\alpha_9,\widetilde{a}_{0 2}^+, \widetilde{b}_{1 0}^+, \widetilde{b}_{1 1}^+), & 
\widetilde{\Psi}_{12} & = \alpha_{9} +  g_{12}(\alpha_7,\alpha_8,\alpha_9,\widetilde{a}_{0 2}^+, \widetilde{b}_{1 0}^+, \widetilde{b}_{1 1}^+), 
\end{aligned}		
\end{equation} 
with $g_i$'s having no linear terms. 
We consider the following blow up of the parameters: 
\begin{equation*}\label{blowupS4}
\alpha_7 = \alpha^2_9\gamma_7, \quad \alpha_8 = \alpha^2_9 \gamma_8, \quad \widetilde{a}^+_{0 2}  =  \alpha_9z_1,\quad \widetilde{b}^+_{1 0} = \alpha_9 z_2, \quad \widetilde{b}^+_{1 1}  = \alpha_9z_3. 
\end{equation*} 
Then the functions \eqref{Psi_1} turn into 
\begin{equation*}\label{Phi_2}
\begin{aligned}
\widetilde{\Psi}_7 & =  \alpha_9^2h_{7}(\alpha_9, \gamma_7,\gamma_8,z_1, z_2, z_3), & \widetilde{\Psi}_8 & = \alpha_9^2  h_{8}(\alpha_9, \gamma_7,\gamma_8 ,z_1, z_2, z_3), \\
 \widetilde{\Psi}_9 & =  \alpha_9^2 h_{9}(\alpha_9, \gamma_7,\gamma_8, z_1, z_2, z_3),  & \widetilde{\Psi}_{10} & =  \alpha_9^2 h_{10}(\alpha_9,\gamma_7,\gamma_8, z_1, z_2, z_3), \\
\widetilde{\Psi}_{11} & =  \alpha_9^2 h_{11}(\alpha_9, \gamma_7,\gamma_8, z_1, z_2, z_3), & \widetilde{\Psi}_{12} & =  \alpha_9 + \alpha_9^2 h_{12}(\alpha_9, \gamma_7,\gamma_8, z_1, z_2, z_3),
\end{aligned}		
\end{equation*}
for suitable analytic functions $h_i$'s in the variables $\alpha_9, \gamma_7,\gamma_8, z_1, z_2, z_3$. 

Before proceeding, we observe that 
$$
\begin{aligned}
h_i = h_i(\alpha_9, \gamma_7,\gamma_8, z_1, z_2, z_3) & = h_{i,0}(\gamma_7,\gamma_8, z_1, z_2, z_3) + \alpha_9 h_{i,1}(\gamma_7,\gamma_8, z_1, z_2, z_3) \\
& \phantom{={}} + \alpha_9^2 h_{i,2}(\gamma_7,\gamma_8, z_1, z_2, z_3) + \cdots, \\
\end{aligned}
$$
for $i = 7, \ldots, 11$, such that the actual expression of $h_{i,0}$ is completely determined by the second order we have calculated above following the algorithm of Section~\ref{section:averaging} for $i = 7, 9, 11$, and for the the first and second order for $i = 8, 10$. 
We do not have the expressions of $h_{i,j}$, for $j \geq 1$, unless we calculate more orders above. 
But anyway for what we want we do not need them. 

It turns out that the quadratic system $h_{7,0} = h_{8,0} = h_{9,0} = h_{10,0} = h_{11,0} = 0$ has a rational solution, given approximately (the exact expressions are too big to write down here) by 
$$
\begin{aligned} 
& \gamma_7^* \approx -1.755\times 10^7, \quad \gamma_8^* \approx -1.318\times10^9,\\
& z_1^* \approx 0.8838, \quad z_2^* \approx 0.09214, \quad z_3^* \approx -0.08745. 
\end{aligned}
$$
Further, the Jacobian determinant of the map $\left(h_{7,0}, h_{8,0}, h_{9,0}, h_{10,0}, h_{11,0}\right): \R^5 \to \R^5$ calculated at this solution $(\gamma_7^*, \gamma_8^*,z_{1}^*,z_{2}^*,z_3^*)$ is different from zero, that is, the intersection of the varieties $h_{7,0} = 0$, $h_{8,0} = 0$, $h_{9,0} = 0$, $h_{10,0} = 0,$ and $h_{11,0} = 0$ is transversal at $(\gamma_7^*, \gamma_8^*,z_{1}^*,z_{2}^*,z_3^*)$. 

Therefore, it follows by the implicit function theorem applied to the equation 
$$
h_7 = h_8 = h_9 = h_{10} = h_{11} = 0, 
$$ 
in the point $(\alpha_9, \gamma_7,\gamma_8, z_1, z_2, z_3) = (0, \gamma_7*, \gamma_8*, z_1^*, z_2^*, z_3^*)$ that we can analytically isolate $\alpha_9 = \alpha_9(\gamma_7,\gamma_8, z_1, z_2, z_3)$ such that this equation remains true along the graphic $\left(\alpha_9(\gamma_7,\gamma_8, z_1, z_2, z_3), \gamma_7, \gamma_8, z_1, z_2, z_3\right)$ for $(\gamma_7, \gamma_8, z_1, z_2, z_3)$ close enough to the point $(\gamma_7^*, \gamma_8^*, z_1^*, z_2^*, z_3^*)$. 
It is clear that $\alpha_9 \neq 0$ if $(\gamma_7, \gamma_8, z_1, z_2, z_3) \neq (\gamma_7^*, \gamma_8^*, z_1^*, z_2^*, z_3^*)$. 
In particular, the intersection $h_7 = h_8 = h_9 = h_{10} = h_{11} = 0$ is transversal along this graphic. 
Moreover, from the expression of $h_{12}$, it follows that by shrinking the neighborhood of $(\gamma_7^*, \gamma_8^*, z_1^*, z_2^*, z_3^*)$ if necessary, we can guarantee that $|h_{12}/\alpha_9| > 1/2$ there.

This means that for any given $\beta_7, \ldots, \beta_{11}$ close enough to zero (but non-zero) we are able to analytically find $\gamma_7,\gamma_8, z_1, z_2, z_3$ in terms of $\beta_7, \ldots, \beta_{11}$ and so $\alpha_9$ depending analytically on them, such that 
$$
h_i = \beta_i,\ i = 7,\ldots, 11, \quad h_{12} \neq 0. 
$$
By considering the problem with the independent variables $\alpha_i$, $i = 1, \ldots, 6$, and assuming we have used the implicit function theorem a priory solving successively $\widetilde{\Psi}_1 = \widetilde{\beta}_1$ in $\alpha_1$ depending on $\widetilde{\beta}_1, \alpha_2, \ldots$, then $\widetilde{\Psi}_2 = \widetilde{\beta}_2$ in $\alpha_2$ depending on the variables $\widetilde{\beta}_1, \widetilde{\beta}_2, \alpha_3, \ldots$, until $\widetilde{\Psi}_6 = \widetilde{\beta}_6$, and then considering the blow up $\widetilde{\beta}_i = \alpha_9^2 \beta_i$, $i = 1,\ldots, 6$, we can finally write 
$$
\begin{aligned}
\Delta_{\beta}(r) & = \alpha_9^2 \beta_1 \big(1 + \CO(r^{12})\big) r + \alpha_9^2 \beta_{2} \big(1 + \CO(r^{11})\big) r^2 + \cdots + \alpha_9^2 \beta_{11} \big(1 + \CO(r^2)\big) r^{11} \\
& \phantom{={}} + \alpha_9 \big( 1 + \CO(\alpha_9) + \CO(r) \big) r^{12}, 
\end{aligned}
$$ 
for free $\beta_1, \ldots, \beta_{11}$ close enough to $0$ (all the other parameters set to zero). 
We are then able to find suitable $\beta_1,  \ldots, \beta_{11}$ and then $\alpha_9$ such that $\Delta_{\beta}(r)$ has at least $11$ zeros. 
That is, for these parameters, the original system has at least $12$ limit cycles. 

\smallskip

Now we will point out the difference in the analysis of the case $S_3$. 
For this case, we have, according to \eqref{F_k}: $f^\pm_0(\theta,r) = 4\cos\theta\left(3\cos^2\theta-4\right)/3$, $g^\pm_0(\theta,r) = -4\sin\theta\cos^2\theta$, 
\begin{equation*}
\begin{aligned}
3 f^\pm_1(\theta,r)&= 4 \Big(b_{2 0}^{\pm} \cos^4 \theta + \big(b_{1 1}^{\pm} - 4 a_{2 0}^{\pm}\big) \cos^3 \theta \sin \theta + \big(b_{0 2}^{\pm} -4 a_{1 1}^{\pm} \big) \cos^2 \theta \sin^2 \theta \\ 
& \phantom{={}} - 4 a_{0 2}^{\pm} \cos \theta \sin^3 \theta \Big) r^3 + \Big(\big(3 a_{2 0}^{\pm} + 4 b_{1 0}^{\pm}\big) \cos^3 \theta + \big(3 a_{1 1}^{\pm} - 16 a_{1 0}^{\pm} + 4 b_{0 1}^{\pm} \\ 
& \phantom{={}} + 3 b_{2 0}^{\pm}\big) \cos^2 \theta \sin \theta + \big(3 a_{0 2}^{\pm} - 16 a_{0 1}^{\pm} + 3 b_{1 1}^{\pm}\big) \cos \theta \sin^2 \theta + 3 b_{0 2} \sin^3 \theta  \Big) r^2 \\ 
&\phantom{={}} + 3 \Big(a_{1 0}^{\pm} \cos^2 \theta + \big(a_{0 1}^{\pm} + b_{1 0}^{\pm}\big) \cos \theta \sin \theta + b_{0 1}^{\pm} \sin^2 \theta \Big) r.
	\end{aligned}
\end{equation*}
Again, due to the size of the expressions of $f_2^{\pm}$, we will not present them explicitly. Following the same idea as in the case of $S_4$ right above, after calculating the expressions of ${\psi}_{i,j}$, $i = 1,2$, $j = 1,\ldots, N$, we eliminate $\varepsilon$ by means of the reparametrization \eqref{elimina} obtaining homogeneus polynomials of degrees $1$ and $2$, respectively, i.e. $\widetilde{\psi}_{1,j}$ and $\widetilde{\psi}_{2, j}$, $j = 1,\ldots, N$. 
Then after solving in sequence 
$$
\widetilde{\psi}_{1,k} = \alpha_k, \ k = 1,\dots,6, \quad \widetilde{\psi}_{1,8} = \alpha_7, \quad \widetilde{\psi}_{1,10} = \alpha_{8}, \quad \widetilde{\psi}_{1,12} = \alpha_{9},
$$
in the variables $\widetilde{a}^-_{10}$, $\widetilde{b}^-_{01}$, $\widetilde{a}^-_{11}$, $\widetilde{a}^-_{0 2}$, $\widetilde{b}^-_{2 0}$, $\widetilde{a}^-_{2 0}$, $\widetilde{a}^+_{1 1}$, $\widetilde{b}^+_{2 0}$, $\widetilde{b}^+_{0 1}$, respectively, it follows that we can rewrite the $\widetilde{\psi}_{1, j}$'s as 
\begin{equation*}\label{linear_o1_S3}
\begin{aligned}
\widetilde{\psi}_{1,k} & =\alpha_k,\ k=1,\ldots, 6,\quad \widetilde{\psi}_{1,8} = \alpha_7, \quad \widetilde{\psi}_{1,10} = \alpha_{8} \quad \widetilde{\psi}_{1,12} = \alpha_{9},\\
\widetilde{\psi}_{1,7} & = -\tfrac{123396623697969152}{2780914306640625} \alpha_1 +\tfrac{47519539134464}{274658203125} \alpha_2 -\tfrac{69377982464}{263671875} \alpha_3 +\tfrac{10399227904}{52734375} \alpha_4 \\
& \phantom{={}} -\tfrac{10598144}{140625} \alpha_5 + \tfrac{5248}{375}\alpha_6,\\
\widetilde{\psi}_{1,9} & = +\tfrac{79291845609546636591104}{14483928680419921875} \alpha_1 -\tfrac{22147997279223453581312}{1042842864990234375} \alpha_2 \\
& \phantom{={}} +\tfrac{88803259464124727296}{2780914306640625} \alpha_3 -\tfrac{1282919781892096}{54931640625} \alpha_4 +\tfrac{55373245251584}{6591796875} \alpha_5 \\
& \phantom{={}} -\tfrac{898318336}{703125} \alpha_6 + \tfrac{1312}{75} \alpha_7 ,\\
\widetilde{\psi}_{1,11} & = -\tfrac{16497681899886282309893372248064}{18331222236156463623046875} \alpha_1 +\tfrac{170572198789950520780428148736 }{48883259296417236328125}\alpha_2 \\
& \phantom{={}} -\tfrac{683024555783581807897739264}{130355358123779296875} \alpha_3 +\tfrac{1328662146680609176551424}{347614288330078125} \alpha_4  \\
& \phantom{={}} -\tfrac{3801923790886678822912}{2780914306640625} \alpha_5 + \tfrac{169088920535957504}{823974609375} \alpha_6 -\tfrac{123758313472 }{52734375}\alpha_7 +\tfrac{2624}{125} \alpha_8, 
\end{aligned}
\end{equation*}
and $\widetilde{\psi}_{1,12}$, $\widetilde{\psi}_{1,13}$, $\widetilde{\psi}_{1,14}$ and $\widetilde{\psi}_{1,15}$ only depending on $\alpha_1, \ldots, \alpha_9$. 

Then acting as above, after a suitable linear change of coordinates, we can eliminate the linear terms from $\widetilde{\Psi}_7$, $\widetilde{\Psi}_{9}$, and $\widetilde{\Psi}_{11}$. 
Further, without loss of generality, we take the condition $\alpha_j=0$, $j = 1,\ldots,6$, obtaining similar equations as \eqref{Psi_1}. 

The difference here is that even making the blow up by using $\alpha_9$, we do not obtain the independence of the intersections of the varieties $h_7 = \cdots = h_{11} = 0$ as for the case $S_4$. 
This means we can not go until $\widetilde{\Psi}_{12}$ and so we invariably will get less limit cycles. 
The result will be even weaker, because we could not reach $\widetilde{\Psi}_{10}$. 
Anyway we proceed in order to illustrate the method. 
We take the simplification 
$$
\widetilde{a}_{0 1}^+ = \widetilde{a}_{1 0}^+ = \widetilde{a}_{0 2}^+ = \widetilde{a}_{2 0}^+ = \widetilde{b}_{0 2}^- = \widetilde{b}_{1 0}^- = \widetilde{b}_{1 1}^- = \widetilde{b}_{0 2}^+ = \widetilde{b}_{1 0}^+ = \widetilde{b}_{1 1}^+ = \alpha_9 = 0
$$ 
and the blow up 
\begin{equation*}\label{blowupS3}
\alpha_7 = \alpha^2_8 \gamma_7, \quad \widetilde{a}^-_{0 1} = \alpha_8z_1,  
\end{equation*} 
so that we can write
\begin{equation*}\label{Phi_5}
\widetilde{\Psi}_7  =  \alpha_8^2h_{7}(\gamma_7,z_1), \quad \widetilde{\Psi}_8  =  \alpha_8^2  h_{8}(\gamma_7,z_1), \quad \widetilde{\Psi}_9  = \alpha_8^2 h_{9}(\gamma_7,z_1), 	
\end{equation*}
for suitable analytic functions $h_7$, $h_8$, and $h_9$, expanded as: 
$$
h_i = h_i(\alpha_8, \gamma_7,z_1) = h_{i,0}(\gamma_7, z_1) + \alpha_8 h_{i,1}(\gamma_7, z_1) + \alpha_8^2 h_{i,2}(\gamma_7,z_1) + \cdots, 
$$
for $i = 7, 8, 9$. The algebraic system $h_{7,0} = h_{8, 0} = 0$ has the rational solution $(\gamma_7^*, z_1^*)$, given approximately by 
$$
(\beta_7^*,z_{1}^*)\approx( 1.403409714\times 10^{12}, -1.862257817\times 10^4). 
$$
This intersection is transversal as the Jacobian determinant of the map $\left(h_{7,0}, h_{8,0} \right)$ with respect to $(\gamma_7,z_{1})$ evaluated at $(\gamma_7^*,z_{1}^*)$ does not vanish. 
Further, we are able to prove that $h_{9,0}(\gamma_7^*, z_1^*)\neq0$. 
Therefore the result follows with the same reasoning as for the case $S_4$. 
\end{proof}

\subsection{Second order for $Z$ with $Z^+ = S^1$ and $Z^- = S^2$, for the straight line $\tau = 1/2$}\label{section:S1S2} 
In this section we prove Theorem~\ref{theorem:main} by pushing forward through order $2$ the calculations initiated in the proof of Theorem~\ref{theorem:main3}, for the fixed straight line given by $\tau = 1/2$. 

\begin{proof}[Proof of Theorem~\ref{theorem:main}] 	
Recall that the origin of the non-perturbed system is a non-degenerate center, so that we can follow the algorithm of Section~\ref{subsection_aproximation} as in the proof of Theorem~\ref{theorem:main3}. 
The expressions of $f^+_0$, $g^+_0$, $f^+_1$ for $S_1$ are given in \eqref{case_S1}, whereas the expressions of $f_0^-$, $g_0^-,$ and $f_1^-$ for $S_2$ are in \eqref{case_S2_fg}. 
The expressions of $f_2^+$ and $f_2^-$ for $S_1$ and $S_2$, respectively, we shall not present explicitly due to the big size of them. 

By fixing $\tau = 1/2$, we directly get from \eqref{ordemum_S1S2} the expressions of $\psi_{1, j}$, $j = 1\ldots, N$. 
Then in order to get the  expressions of $\psi_{2, j}$, $j = 1, \ldots, N$ we follow the calculations of Section~\ref{subsection_aproximation} as in the preceding section: 
\begin{equation*}\label{coef_o222}
\begin{aligned}
\psi_{2,1} & =\pi \Big(  ( b^+_{0 1} + {a^+_{1 0}})^2 \pi - 2 ( b^+_{0 1} + {a^+_{1 0}}) (b_{1 0}^+ - a_{0 1}^+)  \\
& \phantom{={}} - ( {b^-_{0 1}} + {a^-_{1 0}})^2 \pi - 2 ( {b^-_{0 1}} + {a^-_{1 0}}) (b_{1 0}^- - a_{0 1}^-) \Big)/4
\end{aligned}
\end{equation*}
and $\psi_{2, j} = \Omega_{2,j}\left(a^\pm_{\eta}, b^\pm_{\eta}\right)$, 
where $\Omega_{2,j}$ is a suitable homogeneous polynomials of degree $2$ in $a^\pm_{\eta}, b^\pm_{\eta}$, $|\eta|\leq 2$, for $j = 2,\ldots, N$. 

Now, as in the preceding section, we eliminate $\varepsilon$ of $\Delta(r,\varepsilon)$ by redefining the perturbative coefficients according to \eqref{elimina}, obtaining $\Delta(r) = \sum_{j = 1}^{\infty} \widetilde{\Psi}_j r^j$ with $\widetilde{\Psi}_j = \sum_{i = 1}^{\infty} \widetilde{\psi}_{i, j}$, where $\widetilde{\psi}_{i, j}$ is the homogeneous polynomial of degree $i$ in the variables $\widetilde{a}_{\eta}$, $\widetilde{b}_{\eta}$ obtained from $\varepsilon^i\psi_{i,j}$ by the redefinition \eqref{elimina}. 

Then we solve in sequence the linear terms $\widetilde{\psi}_{1,j} = \alpha_j$, $j = 1,\dots,8$, $\widetilde{\psi}_{1,10} = \alpha_9$ and $\widetilde{\psi}_{1,12} = \alpha_{10}$ in $\widetilde{a}^-_{1 0}$, $\widetilde{b}^-_{1 0}$, $\widetilde{a}^-_{11}$, $\widetilde{a}^-_{0 2}$, $\widetilde{b}^+_{2 0}$, $\widetilde{a}^-_{2 0}$, $\widetilde{b}^-_{2 0}$, $\widetilde{a}^+_{0 1}$, $\widetilde{a}^+_{1 0}$, and $\widetilde{a}_{1 1}^+$, respectively, obtaining 
\begin{equation*}\label{linear_o1S1S2}
\begin{aligned}
\widetilde{\psi}_{1,k} & =\alpha_k,\ k=1,\ldots, 8,\quad \widetilde{\psi}_{1,10}=\alpha_9, \quad \widetilde{\psi}_{1,12}=\alpha_{10},\\
\widetilde{\psi}_{1,9} & =\tfrac{27456}{390625}\alpha_1 -\tfrac{36752}{78125}\alpha_2 -\tfrac{22027 }{3125}\alpha_3 -\tfrac{163833}{6250}\alpha_4 \\ &\phantom{={}}-\tfrac{94071}{2000}\alpha_5 -\tfrac{5856 }{125}\alpha_6 -\tfrac{663}{25}\alpha_7  -8 \alpha_8,\\
\widetilde{\psi}_{1,11}& = -\tfrac{26849248}{9765625}\alpha_1 + \tfrac{36771864}{1953125}\alpha_2 +\tfrac{43321833}{156250}\alpha_3 +\tfrac{12752699}{12500}\alpha_4  \\
&\phantom{={}} +\tfrac{900727361 }{500000}\alpha_5+\tfrac{5440604}{3125}\alpha_6 +\tfrac{4604183}{5000}\alpha_7 +\tfrac{5632}{25}\alpha_8 -\tfrac{48}{5}\alpha_9,
\end{aligned}
\end{equation*}
and $\widetilde{\psi}_{1,j}$, $j=12,\ldots, 15$ depending only on $\alpha_j$'s. 

We again pursue a suitable linear change of variables in order to eliminate the linear terms $\widetilde{\psi}_{1, 9}$ and $\widetilde{\psi}_{1, 11}$ of $\widetilde{\Psi}_9$ and $\widetilde{\Psi}_{11}$. 
Then we take without loss of generality $\alpha_j = 0$, $j = 1,\ldots,8$. 
Further, we make the following reduction on the number of variables: 
$$
\widetilde{a}_{0 2}^+ = \widetilde{a}_{2 0}^+ = \widetilde{b}_{0 1}^- = \widetilde{b}_{0 2}^- = \widetilde{b}_{1 1}^- = \widetilde{b}_{0 1}^+ = \widetilde{b}_{0 2}^+ = \widetilde{b}_{1 0}^+ = 0, 
$$
obtaining in particular that
\begin{equation*}\label{Psi_3}
\begin{aligned}
\widetilde{\Psi}_9 & = g_{9}(\alpha_9,\alpha_{10}, \widetilde{a}_{0 1}^-, \widetilde{b}_{1 1}^+),  & \widetilde{\Psi}_{10} & =  \alpha_9 + g_{10}(\alpha_9,\alpha_{10}, \widetilde{a}_{0 1}^-, \widetilde{b}_{1 1}^+),  \\
\widetilde{\Psi}_{11} & = g_{11}(\alpha_9,\alpha_{10}, \widetilde{a}_{0 1}^-, \widetilde{b}_{1 1}^+),  & \widetilde{\Psi}_{12} & =  \alpha_{10} + g_{12}(\alpha_9,\alpha_{10}, \widetilde{a}_{0 1}^-, \widetilde{b}_{1 1}^+), 
\end{aligned}		
\end{equation*}
with $g_{9}$, $g_{10}$, $g_{11},$ and $g_{12}$ suitable analytic functions beginning with order $2$. 

Considering the blow up of the parameters given by
\begin{equation*}\label{blowupS1S2} 
\alpha_9 = \alpha^2_{10} \gamma_9, \quad \widetilde{a}^-_{0 1} = \alpha_{10} z_1, \quad \widetilde{b}^+_{1 1} = \alpha_{10} z_{2}, 
\end{equation*} 
we get 
\begin{equation*}\label{Phi_4}
\begin{aligned}
\widetilde{\Psi}_9 & = \alpha_{10}^2 h_{9}(\gamma_9, \alpha_{10}, z_1, z_2),  &  \widetilde{\Psi}_{10} & = \alpha_{10}^2 h_{10}(\gamma_9, \alpha_{10}, z_1, z_2),  \\
\widetilde{\Psi}_{11} & = \alpha_{10}^2 h_{11}(\gamma_9, \alpha_{10}, z_1, z_2),  & \widetilde{\Psi}_{12} & = \alpha_{10} + \alpha_{10}^2 h_{12}(\gamma_9, \alpha_{10}, z_1, z_2), \\
\end{aligned}		
\end{equation*}
with 
$$
\begin{aligned}
h_i (\gamma_9, \alpha_{10}, z_1, z_2) & = h_{i,0}(\gamma_9, z_1, z_2) + \alpha_{10} h_{i,1}(\gamma_9, z_1, z_2)  + \alpha_{10}^2 h_{i,2}(\gamma_9, z_1, z_2) + \cdots,  
\end{aligned}
$$
for $i = 9, 10, 11, 12$. 

Since the algebraic system 
$$
h_{9,0}(\gamma_0, z_1, z_2) = h_{10, 0}(\gamma_9, z_1, z_2) = h_{11, 0}(\gamma_9, z_1, z_2) = 0 
$$ 
has a rational solution 
$$
(\gamma_{9}^*, z_1^*, z_2^*) \approx( -1.267678465\times 10^{11}, -8.373115792\times10^4, 5.752432052\times 10^4), 
$$ 
and this solution is transversal (because we Jacobian determinant of the map defined by $(h_{9,0}, h_{10,0}, h_{11,0})$ calculated at $(\gamma_{9}^*, z_1^*, z_2^*)$ is nonzero), it follows exactly as in the proof of the preceding section (for $S_4$) that we can guarantee at least $11$ positive zeros of $\Delta_{\beta}(r)$ (observe that $\widetilde{\Psi}_{12} \approx 1$ if $(\gamma_9, \alpha_{10}, z_1, z_2)$ is close to $(\gamma_9^*, 0, z_1^*, z_2^*)$). 
So we get at least $12$ limit cycles, finishing the proof of the theorem. 
\end{proof}

\section{acknowledgements} 
This work has been realized thanks to the Brazilian S\~ao Paulo Research Foundation FAPESP (grants 2019/07316-0, 2020/14498-4, 2021/14987-8, 2022/14484-9 and 2023/00376-2); the Brazilian CAPES Agency (Coordena\c{c}\~ao de Aperfei\c{c}oamento de Pessoal de N\'ivel Superior - Finance Code 001); the Catalan AGAUR Agency (grant 2021 SGR 00113); the Spanish AEI agency (grants PID2019-104658GB-I00 and CEX2020-001084-M), and the European Union's Horizon 2020 research and innovation program (grant Dynamics-H2020-MSCA-RISE-2017-777911). 

\section{Conflict of Interest} 
The authors have no conflicts of interest to declare. 

\bibliographystyle{abbrv}
\bibliography{biblio}
\end{document}